\documentclass[14pt]{extarticle}
\usepackage{extsizes}
\usepackage[a4paper, margin=2cm]{geometry}
\usepackage{subfiles}
\newcommand{\onlyinsubfile}[1]{#1}
\newcommand{\notinsubfile}[1]{}

\everymath{\displaystyle}

\usepackage{amssymb}
\usepackage{amsfonts}
\usepackage{amsthm}
\usepackage{comment}
\usepackage{faktor}
\usepackage{mathtools}
\usepackage{stmaryrd}
\usepackage{thmtools}
\usepackage{tikz-cd}
\usepackage{url}
\usepackage{xifthen}

\usepackage{hyperref}
\usepackage[nameinlink]{cleveref}

\newtheoremstyle{kkkk}
{\topsep}
{0pt}
{\hangindent=2em}
{0pt}
{\bfseries}
{:}
{\newline}
{}

\newtheoremstyle{llll}
{\topsep}
{0pt}
{\hangindent=1em}
{0pt}
{\bfseries}
{}
{\newline}
{}

\newtheoremstyle{pf}
{0pt}
{\topsep}
{\hangindent=2em}
{0pt}
{\bfseries}
{:}
{\newline}
{}

\theoremstyle{kkkk}
\newtheorem{jthm}{Theorem}[section]

\newtheorem{jdef}[jthm]{Definition}

\newtheorem{jlem}[jthm]{Lemma}

\newtheorem{jprop}[jthm]{Proposition}
\newtheorem{jcor}[jthm]{Corollary}

\theoremstyle{pf}
\newtheorem*{jpf}{Proof}

\theoremstyle{llll}
\newtheorem{jexx}[jthm]{Example}

\declaretheoremstyle[
  spaceabove=0pt,
  spacebelow=0pt,
  headfont=\bfseries,
  bodyfont=\hangindent=2em,
  headpunct={},
  postheadspace=\newline,
  notefont=\bfseries,
  notebraces={of~}{:},
]{ofstyle}
\declaretheorem[name={Proof},style=ofstyle, unnumbered]{jpfof}

\newcommand{\Cdim}{\mathrm{Cdim}}
\newcommand{\gr}{\mathrm{gr}\,}

\newcommand{\Ker}{\textrm{Ker\,}}

\newcommand{\overbar}[1]{\mkern 1.5mu\overline{\mkern-1.5mu#1\mkern-1.5mu}\mkern 1.5mu}
\newcommand{\jref}[2][]{\ifthenelse{\isempty{#1}}{[\textbf{\small{#2}}]}{[{\small{\textbf{#2}, #1}}]}}
\newcommand{\jcite}[2][]{\notinsubfile{\ifthenelse{\isempty{#1}}{\cite{#2}}{\cite[#1]{#2}}}\onlyinsubfile{\jref[#1]{#2}}}
\newcommand{\Paskunas}{Pa{\v{s}}k{\={u}}nas }

\begin{document}

\renewcommand{\onlyinsubfile}[1]{}
\renewcommand{\notinsubfile}[1]{#1}

\begin{center}
\LARGE \textbf{The canonical dimension of modules for Iwasawa algebras}\vspace{10pt}\\
\large James Timmins \vspace{-10pt}
\end{center}

\begin{centering} \subsection*{Abstract} \end{centering}

Let $F$ be a non-trivial finite extension of the $p$-adic numbers, and $G$ be a compact $p$-adic Lie group whose Lie algebra is isomorphic to a split semisimple $F$-Lie algebra. We prove that the mod $p$ Iwasawa algebra of $G$ has no modules of canonical dimension one. One consequence is a new upper bound on the Krull dimension of the Iwasawa algebra. We also prove a canonical dimension-theoretic criterion for a mod $p$ smooth admissible representation to be of finite length. Combining our results shows that any smooth admissible representation of $GL_n(F)$, with central character, has finite length if its dual has canonical dimension two.\\

\section*{Introduction}

Let $p$ be a prime number and $k$ be a field of characteristic $p$. When $G$ is a $p$-adic Lie group, a very important class of its representations over $k$ are those known as smooth admissible. The study of smooth admissible representations is crucial to the Langlands programme when $G$ is a general linear group, for example. When the vector spaces $V$ on which the representations are defined have the discrete topology, the smooth representations of $G$ are exactly the continuous representations. The (smooth) admissible representations are those with the property that the subspaces $V^K$ of vectors fixed by an open subgroup $K$ are always finite-dimensional.\\

The vector space dimensions of the invariant subspaces $V^K$ are an important structural property of an admissible representation $V$. Morra computed many of these dimensions for irreducible representations of $GL_2(\mathbb{Q}_p)$ in \jcite{morra2013invariant}. Typically, however, the integers $\dim_k V^K$ are unknown, particularly when $V$ is a supersingular representation of a general linear group.\\

On the other hand, the admissible representations of $G$ are characterised as those whose $k$-linear duals are finitely-generated modules over the Iwasawa algebra of a compact open subgroup, see \Cref{Holonomic representations}. This means the theory of admissible representations is tightly controlled by these Noetherian rings. Because Iwasawa algebras are moreover Auslander-Gorenstein, by \jcite{venjakob2002structure}, the theory of \emph{canonical dimension} gives a natural way to measure and classify their modules. \\

The canonical dimension, as expressed in \Cref{Cdim def}, is in general an integer-valued homological invariant $\Cdim(M)$ for modules $M$ over a noncommutative Auslander-Gorenstein ring $R$. However when $R=kH$ is an Iwasawa algebra, it is fruitful to notice that the canonical dimension is the dimension of an affine variety associated to $M$, known as the characteristic variety, see \Cref{canonical dimension characteristic ideal thm}. We therefore are inclined to view $\Cdim(M)$ as a measure of the ``size'' of $M$, expressed as an integer lying between zero and $\dim H$.\\

In fact, the canonical dimension directly measures ``growth'' in admissible representations. When $V$ is an admissible representation and $M$ is the dual $kH$-module, Emerton and \Paskunas have proved a formula for the asymptotic behaviour of the dimensions $\dim_k V^K$ of invariant subspaces along a chain of compact open subgroups $K$, \jcite{emertonpaskunas2020density}, which crucially involves $\Cdim(M)$. The canonical dimension also appears, although not by name, in conjectures on completed (co)homology \jcite{calegari2012completed} and in the work of Gee and Newton, \jcite{gee2020patching}.\\

More recently, and partially motivated by the works mentioned above, the canonical dimension has been studied by Breuil, Herzig, Hu, Morra, Schraen, and Wang. Their primary interest has been in representations associated with the mod $p$ cohomology of Shimura varieties, which are believed to realise a mod $p$ (local) Langlands correspondence. In \jcite{breuil2023gkdimension} and \jcite{hu2022modpcohomology}, the canonical dimensions associated to such representations of $GL_2(F)$, for $F$ an unramified finite extension of the $p$-adic numbers $\mathbb{Q}_p$ are determined, and found to be equal to the degree of the field extension $[F:\mathbb{Q}_p]$. In \jcite{breuilXconjectures}, the authors make numerous conjectures towards a mod $p$ Langlands correspondence, and prove many of their conjectures due to the results of \jcite{breuil2023gkdimension} and \jcite{hu2022modpcohomology}.
We will take a complementary perspective to that of the above authors: rather than computing the canonical dimension for certain representations, we fix a group and investigate, for a given integer, whether a representation exists with canonical dimension equal to that value. \\

In this article we let $F$ be any finite field extension of $\mathbb{Q}_p$, and consider the Iwasawa algebras of compact $p$-adic Lie groups such as $G=SL_n(\mathcal{O}_F)$, and those corresponding to other semisimple types. Our central result is that when $F \neq \mathbb{Q}_p$, the Iwasawa algebra has no modules of canonical dimension one, \Cref{Main result canonical dimension one thm}.\\ 

From this, we deduce corollaries in two different directions. A restriction on the canonical dimension of modules over a noncommutative ring gives an upper bound on the Krull(-Gabriel-Rentschler) dimension of the ring, \Cref{canonical dimension module range bound prop}. Therefore \Cref{Main result canonical dimension one thm} implies a non-trivial upper bound on the Krull dimension of $kG$, \Cref{Main result Krull dimension bound cor}. This generalises a result of Ardakov, and verifies that a conjecture of Ardakov and Brown holds for the group $SL_2(\mathcal{O}_F)$, when $F$ is a quadratic extension of the $p$-adic numbers.\\

Our second application is to admissible representations of non-compact $p$-adic Lie groups, in particular to representations of the general linear groups $GL_n(F)$. We show that for many split-semisimple $p$-adic Lie groups such as $SL_n(F)$, an admissible representation is of finite length as long as its corresponding canonical dimension takes a certain minimal value, \Cref{Main result minimal-positive implies finite-length thm} and \Cref{Chevalley fd module prop}. One consequence is a characterisation, purely in terms of the canonical dimension, of the admissible finite length representations of $GL_2(\mathbb{Q}_p)$ with central character, see \Cref{Main result GL2(Qp) finite length cor}.\\

In the case $GL_n(F) \neq GL_2(\mathbb{Q}_p)$, we also deduce from \Cref{Main result canonical dimension one thm} that any admissible representation with central character has finite length if its dual has canonical dimension two, \Cref{Main result GLn finite length cor}. Therefore, the admissible representations of $GL_2(\mathbb{Q}_{p^2})$ considered in \jcite{breuil2023gkdimension} and \jcite{hu2022modpcohomology} are necessarily of finite length, as conjectured in \jcite{breuilXconjectures} -- see \Cref{Finite length GLn}.\\

Let us mention that although \Cref{Main result canonical dimension one thm} and \Cref{Main result Krull dimension bound cor} generalise results proved when $F=\mathbb{Q}_p$ in \jcite{ardakov2004krull}, the methods used here are somewhat different. We explicitly rely on having a non-trivial finite field extension of $\mathbb{Q}_p$, and use group structures most clearly found in $SL_2(\mathcal{O}_F)$, which Ardakov's methods cannot apply to.\\

Our methods (like those of \jcite{ardakov2004krull}) do not easily generalise to modules of canonical dimension two or higher. The principal reason for this is that the natural generalisations of the commutative algebra result established in \Cref{commutative algebra} to higher dimensions do not necessarily give sufficient restrictions on the characteristic variety attached to a module. As an illustration of the difficulties, we mention that a module of canonical dimension one has a characteristic variety which corresponds to a projective variety of dimension zero, whereas a module of canonical dimension at least two corresponds to a projective curve or higher-dimensional variety, which have a more complicated structure.\\

\textit{Acknowledgements.}  I would like to heartily thank the anonymous reviewer for their comments on the article, especially for their suggestion of improved proofs for \Cref{Main result minimal-positive implies finite-length thm}, \Cref{commutative algebra prop}, and \Cref{F fd module prop}. I thank Konstantin Ardakov for many valuable discussions. I would like to thank Nicolas Dupré for his comments, especially the observation that the results of this paper hold when $\mathfrak{g}$ is a semisimple Lie algebra. I further thank Mick Gielen, Benjamin Schraen, James Taylor, and Marie-France Vignéras for their helpful comments. This research was financially supported by an EPSRC studentship and EPSRC grant EP/T018844/1.\\

\section{Main results} \label{Main results}

Let $p$ be a prime number and let $k$ be a field of characteristic $p$. Let $F$ be a finite field extension of $\mathbb{Q}_p$, with ring of integers $\mathcal{O}_F$. We maintain this notation throughout the article. \vspace{2pt}\\

In the first part of this article we consider compact $p$-adic Lie groups $G$ whose Lie algebras $\mathcal{L}(G)$ -- see \Cref{Lie algebra def} -- are semisimple and split over $F$. The completed group algebra, known as the Iwasawa algebra $kG$, is a Noetherian Auslander-Gorenstein ring of injective dimension $\dim G$. Therefore any finitely-generated module for $kG$ has a canonical dimension, lying between zero and $\dim G$ -- see \Cref{Cdim def}. We prove that this value cannot be equal to 1 in many cases.

\begin{restatable}{jthm}{canonicaldimensionone} \label{Main result canonical dimension one thm}
Let $\mathfrak{g}$ be a finite-dimensional $F$-Lie algebra which is semisimple and split over $F$. Let $G$ be a compact $p$-adic Lie group with associated Lie algebra $\mathcal{L}(G) \cong \mathfrak{g}$. If $[F:\mathbb{Q}_p] > 1$, then the Iwasawa algebra $kG$ has no finitely-generated modules of canonical dimension one. \\
\end{restatable}

Throughout the article, isomorphisms of Lie algebras are $\mathbb{Q}_p$-Lie algebra isomorphisms. Any $F$-Lie algebra is viewed as a $\mathbb{Q}_p$-Lie algebra via the usual restriction of scalars.\\

\Cref{Main result canonical dimension one thm} has an immediate consequence for the Krull dimension of the Iwasawa algebra.

\begin{restatable}{jcor}{krulldimensionbound} \label{Main result Krull dimension bound cor}
Let $\mathfrak{g}$ be a finite-dimensional $F$-Lie algebra which is semisimple and split over $F$. Let $G$ be a compact $p$-adic Lie group with associated Lie algebra $\mathcal{L}(G) \cong \mathfrak{g}$. If $[F:\mathbb{Q}_p] > 1$, the Krull dimension of $kG$ is at most $\dim G - 1$.\\
\end{restatable}

In \jcite{brumer1966pseudocompact} it was proved that the global dimension of $kG$ is $\dim G$ (or infinity). Thus \Cref{Main result Krull dimension bound cor} shows that the Iwasawa algebra has Krull dimension strictly less than its global dimension.\\

\Cref{Main result canonical dimension one thm} is proved in \Cref{The characteristic ideal}, and \Cref{Main result Krull dimension bound cor} is then deduced in \Cref{Krull dimension}.\\

\Cref{Main result canonical dimension one thm} is a lower bound on the canonical dimension of infinite-dimensional modules for an Iwasawa algebra. This may be viewed as a first step towards a Bernstein's inequality for representations of $p$-adic Lie groups.\\
In the second part of the article, this perspective allows our results to find application to the representation theory of non-compact $p$-adic Lie groups. When $V$ is an admissible representation of a $p$-adic Lie group $G$, the $k$-linear dual $V^\vee = \mathrm{Hom}_k(V,k)$ is a finitely-generated module for the Iwasawa algebra of any compact open subgroup of $G$.\\
We say that an infinite-dimensional representation $V$ is \textit{holonomic} when the canonical dimension $\Cdim(V^\vee)$ is minimal among infinite-dimensional admissible representations of $G$. We prove a result showing that such a representation is often of finite length. 

\begin{restatable}{jthm}{holonomicmodulesfinitelength} \label{Main result minimal-positive implies finite-length thm}
Let $G$ be a $p$-adic Lie group. Suppose $G$ acts trivially on all of its finite-dimensional smooth representations over $k$. Every holonomic smooth admissible representation of $G$ over $k$ is of finite length.\\
\end{restatable}

For example, \Cref{Main result minimal-positive implies finite-length thm} applies to the special linear groups $SL_n(F)$. Due to the existence of non-trivial finite-dimensional smooth representations, the theorem does not apply to general linear groups, or groups of units of $p$-adic division algebras.\\

However, considering representations with central character allows the application of \Cref{Main result minimal-positive implies finite-length thm}. In the case of $GL_2(\mathbb{Q}_p)$, we deduce a characterisation of the finite-length representations with central character.

\begin{restatable}{jcor}{GLQpfinitelength} \label{Main result GL2(Qp) finite length cor}
The finite length smooth admissible representations $V$ of $GL_2(\mathbb{Q}_p)$ with central character are precisely those with $\Cdim(V^\vee) \leq 1$.\\
\end{restatable}

Combining Theorems \ref{Main result canonical dimension one thm} and \ref{Main result minimal-positive implies finite-length thm}, we also deduce the following.

\begin{restatable}{jcor}{GLnfinitelength} \label{Main result GLn finite length cor}
Let $V$ be a smooth admissible representation of $GL_n(F)$ with central character. If $[F:\mathbb{Q}_p] > 1$ or $n >2$, and $\Cdim(V^\vee)=2$, then $V$ is of finite length.\\
\end{restatable}

The proof of \Cref{Main result minimal-positive implies finite-length thm} may be found in \Cref{Holonomic representations}, with \Cref{Main result GL2(Qp) finite length cor} and \Cref{Main result GLn finite length cor} proved in \Cref{Finite length GLn}.\\

\section{Canonical dimension} 

In this section we give two formulations of canonical dimension for modules over Iwasawa algebras.\\

Recall that the Iwasawa algebra of a compact $p$-adic Lie group $G$ is the completed group algebra,
\[kG = \varprojlim_{U \trianglelefteq_o G} k\left[\faktor{G}{U}\right], \]
where the limit is taken over the open normal subgroups of $G$. We summarise the very nice properties of this noncommutative ring.

\begin{jthm} \label{Iwasawa algebra basic properties thm}
The Iwasawa algebra $kG$ is a left and right Noetherian ring, which is Auslander-Gorenstein of injective dimension $\dim G$. When $G$ has no $p$-torsion elements, $kG$ is moreover Auslander-regular.\\
\end{jthm}

In the case that $k=\mathbb{F}_p$, the result follows from Corollary 7.25 of \jcite{ddms2003analyticpro-p}, Theorem J of \jcite{ardakov2007primeness}, and Theorem 3.26 of \jcite{venjakob2002structure}. See also Lemma 6.2 of \jcite{ardakov2012prime}. We give a proof which extends this to any field $k$ of characteristic $p$.\\

We will occasionally refer to augmentation ideals of $kG$, defined below.

\begin{jdef}
Let $N \leq G$ be a closed normal subgroup of a compact $p$-adic Lie group. Its augmentation ideal is the two-sided ideal $\epsilon_G(N) = \Ker(kG \rightarrow k(G/N))$.\\
\end{jdef}

When $G$ is uniform pro-$p$, as in \jcite[Definition 4.1]{ddms2003analyticpro-p}, the augmentation ideal $\epsilon_G(G) = \Ker(kG \rightarrow k)$ is the unique maximal ideal of $kG$, and in this special case we write $J= \epsilon_G(G)$.

\begin{jprop} \label{uniform pro-p filtration prop}
Let $G$ be a uniform pro-$p$ group. The $J$-adic filtration on $kG$ is complete, and the associated graded ring is isomorphic to a commutative polynomial ring in $\dim G$ variables over $k$. \\
\end{jprop}
\begin{jpf}
When $k=\mathbb{F}_p$, section 7.4 of \jcite{ddms2003analyticpro-p} shows that under the $I$-adic filtration, where $I=\mathbb{F}_p[G](G-1)$, the completion $\widehat{\mathbb{F}_p[G]} = \mathbb{F}_pG$, and 
\[ \gr \mathbb{F}_p[G] \cong \gr \mathbb{F}_pG \cong \mathbb{F}_p[X_1, \dots, X_d], \]
where $d$ is the dimension of $G$, \cite[Theorem 7.24]{ddms2003analyticpro-p}. \\
When $k$ is an arbitrary field extension of $\mathbb{F}_p$, we consider the $I_k$-adic filtration on the group algebra $k[G]$, where
\[ I_k = k \otimes_{\mathbb{F}_p} I = k[G](G-1).\]
Clearly $I_k^n = k \otimes_{\mathbb{F}_p} I^n$, thus the associated graded ring is
\[ \gr k[G] = \gr (k\otimes_{\mathbb{F}_p} \mathbb{F}_p[G]) = k\otimes_{\mathbb{F}_p}\gr \mathbb{F}_p[G] \cong k[X_1, \dots, X_d]. \]
Since $G$ is uniform pro-$p$, it is straightforward that $I_k^{p^n} = k[G](G^{p^n}-1)$ with quotient ring $k[G/G^{p^n}]$, where $G^{p^n} = \{g^{p^n} \mid g \in G\}$. The inverse system $(G^{p^n})_{n \geq 0}$ is cofinal in all open normal subgroups of $G$, and therefore
\[ \widehat{k[G]} = \varprojlim_{n \geq 0} k[G]/I_k^n = \varprojlim_{U \trianglelefteq_o G} k[G/U]= kG. \]
It follows that $kG$ has a complete filtration with associated graded
\[ \gr kG = \gr k[G] \cong k[X_1, \dots, X_d]. \]
In particular, $kG$ is Noetherian by \jcite[Proposition 7.27]{ddms2003analyticpro-p}.\\
Now, define $J_n =  \Ker(kG \rightarrow k[G]/I_k^n)$ for all $n \geq 1$. Clearly $J=J_1$, and a direct computation shows that $(J^m + J_n)/J_n = (J_m + J_n)/J_n$ for all $m, n \geq 0$. Therefore the closure of $J^m$ in $kG$ is $J_m$ (in the induced pseudocompact topology).  But $kG$ is Noetherian, hence Corollary 22.4 of \jcite{schneider2011padic} implies all ideals are closed, so $J_m = J^m$. The result follows. $\square$\\
\end{jpf}

\begin{jpfof}[\Cref{Iwasawa algebra basic properties thm}]
Suppose $G$ is uniform pro-$p$. By \Cref{uniform pro-p filtration prop}, the associated graded ring of $kG$ is (left and right) Noetherian, and Auslander-regular (of global dimension $d$). By \jcite[Proposition 7.27]{ddms2003analyticpro-p}, and \jcite[Proposition II.2.2.1, Theorem III.2.2.5]{lvo1996zariskian}, it follows $kG$ is Noetherian and Auslander-regular of dimension $d=\dim G$.\\
Now let $G$ be an arbitrary compact $p$-adic Lie group. Let $G'$ be an open normal uniform pro-$p$ subgroup. By the arguments of section 2.3 of \jcite{ardakov2006survey}, $kG$ is a crossed product with a finite group over the Noetherian ring $kG'$,
\[ kG \cong kG' * (G/G'), \]
hence $kG$ is Noetherian. Since $kG'$ is Auslander-regular, Lemma 5.4 of \jcite{ardakov2007primeness} shows that $kG$ is Auslander-Gorenstein. By Corollary 5.4 of \jcite{ardakov2007primeness}, the injective dimension of $kG$ is equal to the injective dimension of $kG'$, which in turn equals the global dimension by Remark 6.4 of \jcite{venjakob2002structure}. Thus Theorem 4.1 of \jcite{brumer1966pseudocompact} implies $kG$ has injective dimension $\dim G' = \dim G$, and if $G$ has no $p$-torsion elements, that $kG$ is Auslander-regular. $\square$\\
\end{jpfof}

Since $kG$ is Auslander-Gorenstein, Definition 4.5 of \jcite{levasseur1992regulargraded} gives us the definition of canonical dimension below.

\begin{jdef} \label{Cdim def}
Let $M$ be a finitely-generated left $kG$-module. The canonical dimension of $M$ is
$ \Cdim(M) = \dim G - \mathrm{min} \{ j \geq 0 \mid \mathrm{Ext}_{kG}^j(M,kG) \neq 0 \}$. \\
\end{jdef}

Many properties of canonical dimension are proved in section 4 of \jcite{levasseur1992regulargraded}, and a summary can be found in subsection 5.3 of \jcite{ardakov2006survey}.\\

We will frequently work with an alternative characterisation of canonical dimension, understanding the canonical dimension as the dimension of an affine variety known as the \textit{characteristic variety}. The characteristic variety of a module is given by the \textit{characteristic ideal} in the commutative ring $\gr kG$.

\begin{jdef}
Let $G$ be a uniform pro-$p$ group, and recall the $J$-adic filtration on $kG$. Let $M$ be a finitely-generated $kG$-module. The characteristic ideal of $M$ is the radical ideal
\[ J(M) = \sqrt{\mathrm{Ann}(\gr M)} \subseteq \gr kG, \]
where $\mathrm{Ann}(\gr M) = \{ x \in \gr kG \mid x \cdot m=0 \; \forall m \in \gr M \}$ is the annihilator of the graded module $\gr M$.\\
\end{jdef}

Note that $J(M)$ is an ideal of $\gr kG$, whilst the augmentation ideal $J=\epsilon_G(G)$ is a two-sided ideal of $kG$. Recall now the definition of Krull dimension $\mathcal{K}$ from section 6.2 of \jcite{mcconnell2001noncommutative}.

\begin{jthm} \label{canonical dimension characteristic ideal thm}
Let $G$ be a uniform pro-$p$ group, $M$ be a finitely-generated $kG$-module. The canonical dimension of $M$ is given by
\[ \Cdim(M) = \mathcal{K}\Big(\faktor{\gr kG}{J(M)}\Big). \vspace{5pt}\]
\end{jthm}
The proof of \Cref{canonical dimension characteristic ideal thm} is essentially given in section 5.3 of \jcite{ardakov2006survey}, where it is shown that $\Cdim(M) = \mathcal{K}(\gr M)$, which is easily seen to be equal to $\mathcal{K}\left(\gr kG/J(M)\right)$. This integer can be naturally seen as the dimension of the variety associated to $J(M)$, which we do not define. We will frequently use \Cref{canonical dimension characteristic ideal thm} rather than \Cref{Cdim def} when discussing the canonical dimension of a $kG$-module.\\

To finish this section, we prove that the canonical dimension and the Krull dimension are stable when passing to an open normal subgroup. In this article, for $R$ a ring we denote the finitely-generated \emph{non-zero} left $R$-modules by $\mathrm{mod}(R)$.

\begin{jprop} \label{change of group prop}
Let $G_2$ be a compact $p$-adic Lie group and $G_1 \trianglelefteq_o G_2$ be an open normal subgroup. Then,
\[ \{ \Cdim(M) \mid M \in \mathrm{mod}(kG_1)\} = \{ \Cdim(M) \mid M \in \mathrm{mod}(kG_2) \}, \]
and the Krull dimensions $\mathcal{K}(kG_1) = \mathcal{K}(kG_2)$.\\
\end{jprop}
\begin{jpf} 
By the arguments of section 2.3 of \jcite{ardakov2006survey}, we have that $kG_2$ is isomorphic to a crossed product,
\[ kG_2 \cong kG_1 * (G_2/G_1). \]
Thus $kG_2$ is a finitely-generated free $kG_1$-module, and Corollary 6.5.3 of \jcite{mcconnell2001noncommutative} implies that $\mathcal{K}(kG_1) = \mathcal{K}(kG_2)$.\\
Moreover, Lemma 5.4 of \jcite{ardakov2007primeness} implies that for any finitely-generated $kG_2$-module $M$,
\[ \mathrm{Ext}_{kG_2}^n(M,kG_2) \cong \mathrm{Ext}_{kG_1}^n(M,kG_1), \]
and thus $\Cdim_{kG_2}(M) = \Cdim_{kG_1}(M)$. So
\[ \{ \Cdim(M) \mid M \in \mathrm{mod}(kG_1)\} \supseteq \{ \Cdim(M) \mid M \in \mathrm{mod}(kG_2) \}. \]
Conversely, let $M$ be a finitely-generated $kG_1$-module, and let $N = {kG_2 \otimes_{kG_1} M}$. By Theorem 10.74 of \jcite{rotman2009homologicalalgebra},
\[ \mathrm{Ext}_{kG_2}^n(N,kG_2) \cong \mathrm{Ext}_{kG_1}^n(M,kG_2) \cong \bigoplus_{j=1}^m \mathrm{Ext}_{kG_1}^n(M,kG_1), \]
because $kG_2$ is a free $kG_1$-module of rank $m=|G_2/G_1|$. Thus $\mathrm{Ext}_{kG_2}^n(N,kG_2)$ is zero precisely when $\mathrm{Ext}_{kG_1}^n(M,kG_1)$ is zero, from which it follows that $\Cdim_{kG_2}(N) = \Cdim_{kG_1}(M)$. Note also that $N=0$ only if $M=0$. Therefore,
\[ \{ \Cdim(M) \mid M \in \mathrm{mod}(kG_1)\} = \{ \Cdim(M) \mid M \in \mathrm{mod}(kG_2) \}, \]
as required. $\square$\\
\end{jpf}

\section{The Lie algebra} \label{The Lie algebra}

Recall from section 4.5 of \jcite{ddms2003analyticpro-p} that any uniform pro-$p$ group $G$ has an associated $\mathbb{Z}_p$-Lie algebra $L_G$, which has underlying set $G$. We now extend the definition to any compact $p$-adic Lie group, and the scalars to $\mathbb{Q}_p$.

\begin{jdef} \label{Lie algebra def}
Let $G$ be a compact $p$-adic Lie group, and $H \trianglelefteq_o G$ be an open normal uniform pro-$p$ subgroup. The $\mathbb{Q}_p$-Lie algebra $\mathcal{L}(G) = \mathcal{L}(H) = \mathbb{Q}_p \otimes_{\mathbb{Z}_p} L_{H}$ with the obvious $\mathbb{Q}_p$-linear Lie bracket.\\
\end{jdef}

Any compact $p$-adic Lie group has an open normal uniform pro-$p$ subgroup, by Corollary 8.34 of \jcite{ddms2003analyticpro-p}. Since the open subgroups of $G$ all have finite index, it is easy to check that $\mathcal{L}(G)$ does not depend on choice of $H$ up to isomorphism. We now show that the canonical dimension of modules for the Iwasawa algebra, and its Krull dimension, are determined only by the Lie algebra of $G$.

\begin{jcor} \label{Lie algebra equivalence cor}
Let $G_1$, $G_2$ be compact $p$-adic Lie groups with $\mathcal{L}(G_1) \cong \mathcal{L}(G_2)$. Then
\[ \{ \Cdim(M) \mid M \in \mathrm{mod}(kG_1)\} = \{ \Cdim(M) \mid M \in \mathrm{mod}(kG_2) \}, \]
and the Krull dimensions $\mathcal{K}(kG_1) = \mathcal{K}(kG_2)$.\\
\end{jcor}
\begin{jpf} 
Let $K_1 \trianglelefteq G_1$, $K_2 \trianglelefteq G_2$ be open normal uniform pro-$p$ subgroups of $G_1$, $G_2$. Let $\phi$ be an isomorphism of $\mathbb{Q}_p$-Lie algebras,
\[ \phi : \mathcal{L}(G_1) = \mathbb{Q}_p \otimes_{\mathbb{Z}_p} L_{K_1} \rightarrow \mathbb{Q}_p \otimes_{\mathbb{Z}_p} L_{K_2} = \mathcal{L}(G_2). \]
Then, $\phi(L_{K_1})$, $L_{K_2}$ are $\mathbb{Z}_p$-Lie lattices in $\mathcal{L}(G_2)$, so for some $n \geq 0$,
\[ p^n\phi(L_{K_1}) = \phi(p^nL_{K_1}) \leq L_{K_2}. \]
Now, recall from Theorem 9.10 of \jcite{ddms2003analyticpro-p} that the assignment $G \mapsto L_G$ is an equivalence of categories between uniform pro-$p$ groups and powerful $\mathbb{Z}_p$-Lie algebras. $K_1$ is uniform pro-$p$, so $L_{K_1}$ is a powerful Lie algebra, hence $p^nL_{K_1}$ and $\phi(p^nL_{K_1})$ are powerful. Then $\phi(p^nL_{K_1})$ corresponds to a open uniform pro-$p$ subgroup $K \leq G_2$. The $\mathbb{Z}_p$-Lie algebra $p^nL_{K_1}$ corresponds to $K_1^{p^n} = \{x^{p^n} \mid x \in K_1\} \leq G_1$, and thus the restriction of $\phi$ gives us an isomorphism of compact $p$-adic Lie groups $\psi: K_1^{p^n} \rightarrow K$. Now, let
\[ H = \bigcap_{g \in G_2/K} gKg^{-1}, \]
so $H$ is an open normal subgroup of $G_2$. It follows that $H$ is also normal in $K$. So, we have open normal subgroups $K_1^{p^n} \trianglelefteq_o G_1$, $H \trianglelefteq_o K=\psi(K_1^{p^n})$, and $H \trianglelefteq_o G_2$. By \Cref{change of group prop}, it follows that
\begin{align*} \{ \Cdim(M) \mid M \in \mathrm{mod}(kG_1)\} &= \{ \Cdim(M) \mid M \in \mathrm{mod}(kK_1^{p^n})\}\\ &= \{ \Cdim(M) \mid M \in \mathrm{mod}(kK)\}\\ &= \{ \Cdim(M) \mid M \in \mathrm{mod}(kH)\}\\ &= \{ \Cdim(M) \mid M \in \mathrm{mod}(kG_2)\},
\end{align*}
and $\mathcal{K}(kG_1) = \mathcal{K}(kK_1^{p^n}) = \mathcal{K}(kK) = \mathcal{K}(kH) = \mathcal{K}(kG_2)$. $\square$\\   
\end{jpf}

\begin{jexx} \label{abelian example}
Consider the additive group $\mathcal{O}_F$, and let $d = [F:\mathbb{Q}_p]$. Then the Iwasawa algebra can be identified as $k\mathcal{O}_F \cong k[[T_1, \dots, T_d]]$. Consider the $k\mathcal{O}_F$-modules $M_n = k\mathcal{O}_F/I_n$, where $I_n$ is the ideal generated by the first $n$ elements of $\{T_1, \dots, T_d\}$. It is a straightforward deduction from \Cref{canonical dimension characteristic ideal thm} that $\Cdim(M_n)=d-n$, hence
\[ \{ \Cdim(M) \mid M \in \mathrm{mod}(k\mathcal{O}_F)\} = \{0,1, \dots, d\}. \]
Also, the Krull dimension $\mathcal{K}(k\mathcal{O}_F) = \mathcal{K}(k[[T_1, \dots, T_d]]) = d$.\\
Now, $\mathcal{O}_F$ is an abelian group, hence $\mathcal{L}(\mathcal{O}_F)$ must be the abelian $\mathbb{Q}_p$-Lie algebra of dimension $[F:\mathbb{Q}_p]$. By \Cref{Lie algebra equivalence cor}, the above results apply to any abelian compact $p$-adic Lie group of dimension $[F:\mathbb{Q}_p]$. In particular, 
\[ \{ \Cdim(M) \mid M \in \mathrm{mod}(k\mathcal{O}_F^\times)\} = \{0,1, \dots, d\}, \]
and $\mathcal{K}(k\mathcal{O}_F^\times) = [F:\mathbb{Q}_p]$.\\
\end{jexx}

Motivated by \Cref{Lie algebra equivalence cor}, we now associate a specific group to each split-semisimple $F$-Lie algebra, which will assist us in our later proofs.\\

Throughout this article, $\mathfrak{g}$ will denote a finite-dimensional $F$-Lie algebra which is semisimple and split over $F$.\\

We fix a uniform pro-$p$ group $G_\mathfrak{g}$ with Lie algebra $\mathcal{L}(G_\mathfrak{g}) \cong \mathfrak{g}$, as $\mathbb{Q}_p$-Lie algebras. \Cref{Lie algebra equivalence cor} will then show that proving \Cref{Main result canonical dimension one thm} and \Cref{Main result Krull dimension bound cor} for $G_\mathfrak{g}$ implies these results for arbitrary $G$. We make the choice of $G_{\mathfrak{g}}$ in the following way -- see also section 3.4 of \jcite{ardakov2004krull} for a similar construction.\\

Because $\mathfrak{g}$ is split semisimple over $F$, let $\mathfrak{g} = X_F$, where $X$ is a root system. For any commutative ring $R$, there is a Chevalley basis for the Lie algebra $X_R$, namely an $R$-basis $\mathcal{B}$ such that $\mathbb{Z}\mathcal{B}$ is closed under the Lie bracket. Such a basis was constructed by Chevalley for simple $\mathfrak{g}$, and the semisimple case follows by taking a union. We can then define the adjoint Chevalley group $X(R) \leq \mathrm{Aut}(X_R)$ in the usual way.\\
Consider the $\mathcal{O}_F$-Lie algebra $X_{\mathcal{O}_F}$, and notice that the subalgebra $pX_{\mathcal{O}_F}$ is a powerful $\mathcal{O}_F$-Lie algebra. Let $Y$ be the uniform pro-$p$ group constructed from $pX_{\mathcal{O}_F}$ using the Campbell-Hausdorff formula. Let $\mathrm{Ad}: Y \rightarrow GL(pX_{\mathcal{O}_F})$ be the group homomorphism given by $\mathrm{Ad}(g)(u) = gug^{-1}$.

\begin{jdef} \label{specific G def}
The compact $p$-adic Lie group $G_{\mathfrak{g}} = \mathrm{Ad}(Y)$ is the image of the group homomorphism $\mathrm{Ad}$.\\
\end{jdef}

It is straightforward to show that $\Ker \mathrm{Ad} = Z(Y) = \{e\}$, so $\mathrm{Ad}$ is an injective homomorphism. It is also clear that $G=G_\mathfrak{g}$ has Lie algebra 
\[ \mathcal{L}(G) = \mathbb{Q}_p L_G \cong \mathbb{Q}_p L_Y = \mathbb{Q}_p pX_{\mathcal{O}_F} = X_F = \mathfrak{g}. \]
We define some root subgroups of $G_{\mathfrak{g}}$, corresponding to the root system of the Lie algebra $\mathfrak{g}$.

\begin{jdef} \label{root subgroup def}
The root subgroup corresponding to $\alpha \in X$ is $A_\alpha = \mathrm{Ad}(p(X_{\mathcal{O}_F})_\alpha) \leq G_{\mathfrak{g}}$.\\
\end{jdef}

The root subgroups $A_{\alpha}$ will be used over the following three sections, in the proof of \Cref{Main result canonical dimension one thm}.

\begin{jlem} \label{Chevalley normal lemma}
$G_\mathfrak{g}$ is a normal subgroup of the Chevalley group $X(\mathcal{O}_F)$. Moreover, Ad is $X(O_F)$-equivariant with respect to the defining action on $X_{O_F}$ and conjugation action on $G$.
\end{jlem}
\begin{jpf} 
This is similar to the proof of Lemma 3.13 of \jcite{ardakov2004krull}. Let $\mathcal{B}$ be the Chevalley basis of $X_{\mathcal{O}_F}$. It is straightforward to show that $\mathrm{Ad}(pu\mathcal{O}_F) \subseteq X(\mathcal{O}_F)$ for all $u \in \mathcal{B}$. There exists a $\mathbb{Z}_p$-basis $\{pv_1, \dots, pv_N \}$ of $pX_{\mathcal{O}_F}$ where each element is contained in some $pu\mathcal{O}_F$. Now, $Y$ is generated by the subgroups $pv_j\mathbb{Z}_p$, by Proposition 3.7 and Theorem 9.8 of \jcite{ddms2003analyticpro-p}, and $\mathrm{Ad}(pv_j\mathbb{Z}_p) \subseteq X(\mathcal{O}_F)$ for each $j$. So $G_\mathfrak{g} = \mathrm{Ad}(Y) \subseteq X(\mathcal{O}_F)$.\\
It is then straightforward to check that $G_\mathfrak{g}$ is normal in $X(\mathcal{O}_F)$, using the Steinberg relations.\\
In particular $X(\mathcal{O}_F)$ acts on $G$ by conjugation. Now Exercise 9.10 of \jcite{ddms2003analyticpro-p}, implies $\mathrm{Ad}(pv) = \exp(\mathrm{ad}(pv))$. Let $x \in X(\mathcal{O}_F)$, $v \in X_{\mathcal{O}_F}$, and $u \in pX_{\mathcal{O}_F}$. Now, because the Chevalley group consists of Lie algebra automorphisms, we have $x \cdot [pv, x^{-1}\cdot u] = [p(x \cdot v), u]$, and iterating, we have
\[x \cdot \mathrm{ad}(pv)^n(u) = \mathrm{ad}(p(x \cdot v))^n(u).\]
It follows from the formula for $\exp$ that
\[ x\mathrm{Ad}(pv)x^{-1}(u) = x \cdot \exp(\mathrm{ad}(pv))(x^{-1} \cdot u) 
= \exp(\mathrm{ad}(p(x \cdot v)))(u). \]
Namely,
\[ x\mathrm{Ad}(pv)x^{-1} = \mathrm{Ad}(p(x \cdot v)), \]
as required. $\square$\\
\end{jpf}

\section{The microlocalised module}

Recall the theory of algebraic microlocalisation from Chapter IV of \jcite{lvo1996zariskian}. If $\bar{S} \subseteq \gr kG$ is a multiplicatively closed subset of homogeneous elements, we can extend $kG$ to a microlocalised algebra $Q_{\bar{S}}(kG)$. Similarly, for any finitely-generated $kG$-module $M$, we can extend $M$ to the microlocalisation $Q_{\bar{S}}(M)$.\\ 
We prove that a certain microlocalisation of $M$ is isomorphic to a finite-dimensional vector space.

\begin{jprop} \label{microlocalisation fd prop}
Let $G$ be a uniform pro-$p$ group, and $M$ be a finitely-generated left $kG$-module. Let $g \in G$ be a non-trivial element, let $H = \langle g \rangle \leq G$, and let $\bar{S} = \{ {\sigma(g-1)}^n \mid n \geq 0 \}$. Suppose that $\bar{S}^{-1}\Big(\frac{\gr kG}{J(M)}\Big)$ is a finitely-generated $\bar{S}^{-1}k[\sigma(g-1)]$-module.\\ 
Then the microlocalisation $Q_{\bar{S}}(M)$ is a finite-dimensional $K$-vector space, where $K$ is the field $\mathrm{Frac}(kH) \cong \mathrm{Frac}(k[[g-1]]) = k((g-1))$.\\
\end{jprop}
\begin{jpf} 
Because $J(M)$ is the radical of $\mathrm{Ann}(\gr M)$, and $\gr kG$ is a commutative Noetherian ring, there exists $r \in \mathbb{N}$ such that $J(M)^r$ annihilates $\gr M$. Because $\gr M$ is finitely-generated over $\gr kG$, it follows that there is a surjective homomorphism
\[ \Big(\frac{\gr kG}{J(M)^r} \Big)^n \rightarrow \gr M, \]
hence we have a surjective homomorphism
\[ \bar{S}^{-1}\Big( \frac{\gr kG}{J(M)^r} \Big)^n \rightarrow \bar{S}^{-1}\gr M, \] for some $n \in \mathbb{N}$.
Now, the module $\bar{S}^{-1} \big({\gr kG}/{J(M)^r}\big)$ has a series of subquotients
\[ \bar{S}^{-1} \Big(\frac{\gr kG}{J(M)}\Big), \,\bar{S}^{-1} \Big(\frac{J(M)}{J(M)^2}\Big), \dots, \,\bar{S}^{-1} \Big(\frac{J(M)^{r-1}}{J(M)^r}\Big). \]
Each of these subquotients is finitely-generated over $\bar{S}^{-1}\gr kG$ and annihilated by $\bar{S}^{-1}J(M)$. Therefore each is a finitely-generated module for the ring $\bar{S}^{-1}\big({\gr kG}/{J(M)}\big)$. Because $\bar{S}^{-1}\big({\gr kG}/{J(M)}\big)$ is finitely-generated as a $\bar{S}^{-1}k[\sigma(g-1)]$-module, it follows that each subquotient, and hence $\bar{S}^{-1} \big({\gr kG}/{J(M)^r}\big)$ itself, is finitely-generated over $\bar{S}^{-1}k[\sigma(g-1)]$. Considering the surjective homomorphism above, it follows that $\bar{S}^{-1}\gr M$ is a finitely-generated $\bar{S}^{-1}k[\sigma(g-1)]$-module.\\
Now, the closed subgroup $H$ is isomorphic to $\mathbb{Z}_p$, with $kH \cong k[[g-1]]$. By Proposition 4.1.8 of \jcite{lvo1996zariskian},
\[ \bar{S}^{-1}k[\sigma(g-1)] = \bar{S}^{-1}\gr kH = \gr Q_{\bar{S}}(kH), \]
and we have shown $\bar{S}^{-1}\gr M \cong \gr Q_{\bar{S}}(M)$ is a finitely-generated $\gr Q_{\bar{S}}(kH)$-module. By Proposition 4.1.7 and Theorem 1.5.7 of \jcite{lvo1996zariskian}, the finite-generation can be lifted since $Q_{\bar{S}}(M)$ has a complete filtration, thus $Q_{\bar{S}}(M)$ is a finitely-generated $Q_{\bar{S}}(kH)$-module. But
\[ Q_{\bar{S}}(kH) = \textrm{Frac}(kH) = K \]
is a field, so the statement follows. $\square$\\
\end{jpf}

We combine this with the following commutative algebra proposition, which is proved in \Cref{commutative algebra}.

\begin{jprop} \label{characteristic ideal fg element prop}
Let $G$ be a uniform pro-$p$ group, $M$ be a finitely-generated left $kG$-module of canonical dimension 1. Let $g \in G$ be such that $\sigma(g-1) \not \in J(M)$. Then, $\bar{S}^{-1}\Big(\frac{\gr kG}{J(M)}\Big)$ is a finitely-generated $\bar{S}^{-1}k[\sigma(g-1)]$-module.\\
\end{jprop}

Recall the uniform pro-$p$ group $G_\mathfrak{g}$ and its root subgroups $A_{\alpha}$ from Definitions \ref{specific G def} and \ref{root subgroup def}.

\begin{jcor} \label{non-zero annihilator cor}
Suppose $[F:\mathbb{Q}_p] > 1$. Let $M$ be a finitely-generated left $kG_\mathfrak{g}$-module of canonical dimension 1. Let $\alpha \in X$ be a root. Suppose there is a $g \in A_{\alpha}$ such that $\sigma(g-1) \not \in J(M)$. Then the annihilator $\mathrm{Ann}_{kA_{\alpha}}(Q_{\bar{S}}(M))$ is non-zero.\\
\end{jcor}
\begin{jpf} 
By Propositions \ref{microlocalisation fd prop} and \ref{characteristic ideal fg element prop}, $Q_{\bar{S}}(M)$ is a finite-dimensional vector space over the field $K = k((g-1))$. Now, the root subgroup corresponding to $\alpha$ is an abelian uniform pro-$p$ group,
\[A_{\alpha} \cong \mathcal{O}_F \cong \mathbb{Z}_p^{[F:\mathbb{Q}_p]}. \]
Since $[F:\mathbb{Q}_p] > 1$, let $h \in A_{\alpha}$ be non-trivial such that the closed subgroup of $A_{\alpha}$ generated by $g,h$ is 2-dimensional. Because $g,h$ commute, $h$ acts as a $K$-linear map on $Q_{\bar{S}}(M)$. Since $Q_{\bar{S}}(M)$ is finite-dimensional over $K$, the characteristic polynomial gives a non-zero element of $K[h]$ that annihilates $Q_{\bar{S}}(M)$. By clearing denominators, there exists a non-zero element of $k[[g-1]][h] \leq kA_{\alpha}$ annihilating $Q_{\bar{S}}(M)$. Thus the annihilator ideal is non-zero as claimed. $\square$\\
\end{jpf}

\section{$\Gamma$-invariant ideals}

In this section we recall and slightly extend a result of Ardakov, which in our situation classifies the prime ideals of the Iwasawa algebra $kA_\alpha$ of a root subgroup invariant under the action of a torus.
Recall the following definition from subsection 8.1 of \jcite{ardakov2012prime}.

\begin{jdef} 
Let $A$ be a free abelian pro-$p$ group of finite rank, and $\Gamma$ be a closed subgroup of the continuous automorphisms of $A$. We say that $\Gamma$ acts rationally irreducibly on $A$ if every non-trivial $\Gamma$-invariant closed subgroup of $A$ is open.\\
\end{jdef}

Equivalently, $\Gamma$ acts rationally irreducibly if and only if the $\mathbb{Q}_p$-Lie algebra $\mathcal{L}(A)$ is an irreducible module for the Lie algebra $\mathcal{L}(\Gamma)$. Ardakov has proved the following restriction on $\Gamma$-invariant prime ideals in this situation, see Corollary 8.1 of \jcite{ardakov2012prime}.

\begin{jthm} \label{gamma-invariant prime ideals thm}
Suppose $\Gamma$ acts rationally irreducibly on $A$ and 
\[ (h \cdot g) g^{-1} \in A^p = \{x^p \mid x \in A \} \]
for all $h \in \Gamma, g \in A$. The only non-zero $\Gamma$-invariant prime ideal of $kA$ is the unique maximal ideal.\\
\end{jthm}

We extend this result slightly.

\begin{jprop} \label{gamma-invariant ideals prop}
Suppose $\Gamma$ acts rationally irreducibly on $A$ and $(h \cdot g) g^{-1} \in A^p$ for all $h \in \Gamma, g \in A$. If $I$ is a non-zero $\Gamma$-invariant ideal of $kA$, then $\sqrt{I}$ is maximal. \\
\end{jprop}
\begin{jpf} 
Note that $\Gamma$ acts continuously on $kA$ by ring automorphisms. Let $P$ be a minimal prime over $I$, and $\Gamma_P = \{h \in \Gamma \mid h \cdot P = P \}$. Since $kA$ is Noetherian, let $P$ have generators $x_1, \dots, x_n \in kA$. Then
\[ \Gamma_P = \bigcap_{i=1}^n \phi_{i}^{-1}(P) \]
where $\phi_i$ is the restriction of the action map to $\Gamma \cong \Gamma \times \{x_i\} \rightarrow kA$. Because $P$ is closed in $kA$, each $\phi_{i}^{-1}(P)$ is closed. Hence $\Gamma_P$ is a closed subgroup of $\Gamma$. Now, $h \cdot P$ is a minimal prime over $h \cdot I = I$ for any $h \in \Gamma$, and there are finitely many minimal primes over $I$, because $kA$ is Noetherian. Therefore $\Gamma_P$ must be a finite index closed subgroup of $\Gamma$, so is an open subgroup. Hence the inclusion map $\Gamma_P \rightarrow \Gamma$ induces an isomorphism $\mathcal{L}(\Gamma_P)\cong \mathcal{L}(\Gamma)$, so $\Gamma_P$ acts rationally irreducibly on $A$. Moreover, $(h \cdot g) g^{-1} \in A^p$ for all $h \in \Gamma_P$ because $\Gamma_P \leq \Gamma$.\\
By definition of $\Gamma_P$, $P$ is a non-zero $\Gamma_P$-invariant prime ideal of $kA$. Thus by \Cref{gamma-invariant prime ideals thm}, $P$ is the unique maximal ideal of $kA$. So, the only minimal prime over $I$ is the unique maximal ideal, and hence the radical $\sqrt{I}$ is maximal. $\square$\\
\end{jpf}

We now apply this in the case that $A$ and $\Gamma$ are both subgroups of $G_\mathfrak{g}$, with $\Gamma$ acting by conjugation. Recall $G_\mathfrak{g}$ is uniform pro-$p$ and has $\mathbb{Z}_p$-Lie algebra $pX_{\mathcal{O}_F}$, which is a powerful $\mathcal{O}_F$-Lie algebra.\\

For the remainder of this article, fix $\Gamma$ as follows.

\begin{jdef} 
Let $p\mathfrak{t}' \leq pX_{\mathcal{O}_F}$ be the maximal torus of $pX_{\mathcal{O}_F}$ which corresponds to the root system $X$. The subgroup $\Gamma = \mathrm{Ad}(p\mathfrak{t}') \leq G_\mathfrak{g}$.\\
\end{jdef}

Then $\Gamma$ is an abelian closed subgroup of $G_\mathfrak{g}$. Since $p\mathfrak{t}'$ normalises $p({X_{\mathcal{O}_F}})_\alpha$ for any root $\alpha \in X$, it follows that $\Gamma$ normalises any root subgroup $A_\alpha$.

\begin{jlem} \label{rationally irreducibly lemma}
The group $\Gamma$ acts rationally irreducibly, by conjugation, on $A_{\alpha}$ for any root $\alpha \in X$.\\
\end{jlem}
\begin{jpf} 
We have that $\mathcal{L}(A_{\alpha}) = \mathfrak{g}_{\alpha}$ and $\mathcal{L}(\Gamma) = \mathfrak{t}$ is the maximal torus of $\mathfrak{g}$ corresponding to the root system $X$. Both can be considered as $F$-Lie algebras, with $\dim_F \mathfrak{g}_\alpha = 1$, and $\mathfrak{t}$ acts via scalars given by $\alpha: \mathfrak{t} \rightarrow F$. Therefore $\mathcal{L}(A_\alpha)$ is an irreducible $\mathcal{L}(\Gamma)$-module. So $\Gamma$ acts rationally irreducibly. $\square$ \\
\end{jpf}

\begin{jlem} \label{conjugation calculation lemma}
For all $h \in \Gamma$, $g \in A_{\alpha}$, $hgh^{-1}g^{-1} \in A_{\alpha}^p$, for all roots $\alpha \in X$.\\
\end{jlem}
\begin{jpf} 
This follows from the fact that $L_G = pX_{\mathcal{O}_F}$ is a powerful Lie algebra and that $\Gamma$ normalises $A_{\alpha}$. $\square$ \\
\end{jpf}

Combining Lemmas \ref{rationally irreducibly lemma} and \ref{conjugation calculation lemma} together with \Cref{gamma-invariant ideals prop}, we have the following.

\begin{jcor} \label{gamma-invariant ideal specific cor} 
Let $\alpha \in X$. If $I$ is a non-zero $\Gamma$-invariant ideal of $kA_{\alpha}$, then the radical of $I$ is the augmentation ideal $\sqrt{I} = \epsilon_{A_\alpha}(A_\alpha)$.\\
\end{jcor}

\section{The characteristic ideal} \label{The characteristic ideal}

In this section we prove our first main result, \Cref{Main result canonical dimension one thm}. We do this by finding many elements which must lie in the characteristic ideal of a module of canonical dimension one.

\begin{jthm} \label{elements of characteristic ideal thm}
Let $M$ be a finitely-generated left $kG_\mathfrak{g}$-module of canonical dimension 1. If $[F:\mathbb{Q}_p] > 1$, then for any root $\alpha \in X$ and any $g \in A_{\alpha}$, $\sigma(g-1) \in J(M)$. \\
\end{jthm}
\begin{jpf} 
For a contradiction, suppose not. Then by \Cref{non-zero annihilator cor}, there is an $\alpha \in X$ and a non-trivial $g \in A_{\alpha}$ such that the annihilator ideal
\[ I = \mathrm{Ann}_{kA_{\alpha}}(Q_{\bar{S}}(M)) \trianglelefteq kA_{\alpha} \]
is non-zero, where $\bar{S} = \{\sigma(g-1)^n \mid n \geq 0 \}$. Now, because the microlocalisation $Q_{\bar{S}}(M)$ has an action of $kG$, for any $x \in I$, $h \in \Gamma$, and $m \in Q_{\bar{S}}(M)$, we have
\[ hxh^{-1} \cdot m = hx\cdot(h^{-1}m) = h\cdot 0 = 0, \]
therefore $hxh^{-1} \in I$. So, $I$ is a $\Gamma$-invariant ideal of $kA_{\alpha}$. By \Cref{gamma-invariant ideal specific cor}, $\sqrt{I} = \epsilon_{A_{\alpha}}(A_{\alpha})$ contains $g-1$, so $(g-1)^n \in I$ for some $n \geq 1$. But $Q_{\bar{S}}(M)$ is a vector space over $K = k((g-1))$, by \Cref{microlocalisation fd prop}, and is annihilated by a non-zero element of $K$, thus $Q_{\bar{S}}(M)=0$. So $\bar{S}^{-1} \gr M = \gr Q_{\bar{S}}(M) = 0$, hence $\bar{S} \cap \mathrm{Ann}(\gr M) \neq \emptyset$. It follows that ${\sigma(g-1) \in J(M)}$, a contradiction. $\square$ \\
\end{jpf}

\begin{jcor} \label{automorphism characteristic ideal cor}
Let $M$ be a finitely-generated left $kG_\mathfrak{g}$-module of canonical dimension 1. If $[F:\mathbb{Q}_p] > 1$, then for all continuous group automorphisms $\phi: G_\mathfrak{g} \rightarrow G_\mathfrak{g}$, all roots $\alpha \in X$ and all $g \in A_{\alpha}$, the element $\sigma(\phi(g)-1) \in J(M)$.\\
\end{jcor}
\begin{jpf} 
For ease we write $G=G_\mathfrak{g}$. Any continuous group automorphism $\phi$ extends to a filtered ring automorphism $\widetilde{\phi}: kG \rightarrow kG$, which induces a graded ring automorphism $\bar{\phi}: \gr kG \rightarrow \gr kG$. We then define a $kG$-module $M'$ on the set $M$ via $x \cdot m = \widetilde{\phi}(x)m$. It is easy to see that $J(M') = \bar{\phi}^{-1}(J(M))$, and therefore $M'$ is a module of canonical dimension 1. By \Cref{elements of characteristic ideal thm}, $\sigma(g-1) \in J(M')$ for any $\alpha \in X$, $g \in A_{\alpha}$. Hence $\bar{\phi}(\sigma(g-1)) = \sigma(\phi(g)-1) \in J(M)$. $\square$\\
\end{jpf}

We show that this forces the characteristic ideal $J(M)$ to be too large, by showing that $J(M)$ must contain the entirety of the first graded piece $J/J^2$ of $\gr kG_\mathfrak{g}$.

\begin{jlem} \label{identification lemma}
Let $G=G_\mathfrak{g}$. There is an isomorphism of $k[\mathrm{Aut}(G)]$-modules,
\[ \faktor{J}{J^2} \cong k\otimes_{\mathbb{F}_p} \left( \faktor{G}{G^p}\right) \cong k\otimes_{\mathbb{F}_p}\left( \faktor{pX_{\mathcal{O}_F}}{p^2X_{\mathcal{O}_F}} \right), \]
where $\sigma(g-1) +J^2\mapsto gG^p$ and $\mathrm{Ad}(pv)G^p \mapsto pv + p^2X_{\mathcal{O}_F}$ for $g \in G$, $v \in X_{\mathcal{O}_F}$.\\
\end{jlem}
\begin{jpf} 
It suffices to prove the statement when $k=\mathbb{F}_p$, by then applying the appropriate tensor product.\\
In this case, the first isomorphism is given by Lemma 3.11 of \jcite{ardakov2004krull}. For the second isomorphism, Corollary 4.18 of \jcite{ddms2003analyticpro-p} implies that the identity map $\phi: G \rightarrow L_G$ is an isomorphism of $\mathbb{Z}_p[\mathrm{Aut}(G)]$-modules. The image $\phi(G^p) = pL_G$, and therefore 
\[ \overbar{\phi}:\faktor{G}{G^p} \rightarrow \faktor{L_G}{pL_G} \]
is an isomorphism of $\mathbb{F}_p[\mathrm{Aut}(G)]$-modules. Since $L_G =pX_{\mathcal{O}_F}$, the statement follows. $\square$\\
\end{jpf}

We can now prove \Cref{Main result canonical dimension one thm}.

\notinsubfile{\canonicaldimensionone*}
\onlyinsubfile{
\begin{jthm} \label{Main result canonical dimension one thm}
Let $\mathfrak{g}$ be a finite-dimensional $F$-Lie algebra which is semisimple and split over $F$. Let $G$ be a compact $p$-adic Lie group with associated Lie algebra $\mathcal{L}(G) \cong \mathfrak{g}$. If $[F:\mathbb{Q}_p] > 1$, then the Iwasawa algebra $kG$ has no finitely-generated modules of canonical dimension one. \\
\end{jthm}
}
\begin{jpf} 
By \Cref{Lie algebra equivalence cor}, we may assume that $G = G_\mathfrak{g}$ is the uniform pro-$p$ group of \Cref{specific G def}.\\ Suppose there exists a finitely-generated $kG$-module $M$ with canonical dimension 1. We obtain a contradiction by showing that $J(M)$ contains $J/J^2$, as in this case $J(M)$ is the maximal graded ideal of $\gr kG$, which forces $\Cdim(M)=0$.\vspace{2pt}\\
By \Cref{identification lemma}, let $L$ be the image in $k\otimes_{\mathbb{F}_p}\left( pX_{\mathcal{O}_F}/p^2X_{\mathcal{O}_F} \right)$ of the degree 1 elements of $J(M)$. We adopt the notation of section 3.4 of \jcite{ardakov2004krull}, which naturally extends to semisimple Lie algebras, so that $X_{\mathcal{O}_F}$ is spanned by $\{e_{\alpha}, h_{\alpha} \mid \alpha \in X \}$ and the Chevalley group $X(\mathcal{O}_F)$ is generated by elements $\{ x_\alpha(t) \mid \alpha \in X, t \in \mathcal{O}_F\}$. \vspace{2pt}\\
Since $L$ is a $k$-vector space, it is enough to show that $L$ contains the $\mathcal{O}_F/p\mathcal{O}_F$-Lie algebra $pX_{\mathcal{O}_F}/p^2X_{\mathcal{O}_F}$. By considering the group elements $\mathrm{Ad}(p\lambda e_\alpha) \in A_{\alpha}$, \Cref{identification lemma} and \Cref{elements of characteristic ideal thm} imply that for any $\lambda \in \mathcal{O}_F$ and any $\alpha \in X$,
\[p\lambda e_{\alpha} + p^2X_{\mathcal{O}_F} \in L. \]
By the first part of \Cref{Chevalley normal lemma}, $X(\mathcal{O}_F)$ acts by conjugation automorphisms on $G$, hence on the $k[\mathrm{Aut}(G)]$ modules of \Cref{identification lemma}, and on $L$. Thus by \Cref{automorphism characteristic ideal cor},
\[ x \cdot (p\lambda e_{\alpha} + p^2X_{\mathcal{O}_F}) \in L, \]
for any $x \in X(\mathcal{O}_F)$, $\lambda \in \mathcal{O}_F$ and $\alpha \in X$. By the second part of \Cref{Chevalley normal lemma}, the identification of \Cref{identification lemma} implies $X(\mathcal{O}_F)$ acts on $pX_{\mathcal{O}_F}/p^2X_{\mathcal{O}_F}$ via its standard action on $X_{\mathcal{O}_F}$. Recall that $x_{-\alpha}(1) \in X(\mathcal{O}_F)$ acts on $e_{\alpha} \in X_{\mathcal{O}_F}$ by
\[ x_{-\alpha}(1) \cdot e_{\alpha} = e_{\alpha} + h_{-\alpha} - e_{-\alpha}. \]
Since the images of $p\lambda e_{\alpha}$, $p \lambda e_{-\alpha}$ lie in $L$, it follows that
\[ p \lambda h_{-\alpha} + p^2X_{\mathcal{O}_F} = x_{-\alpha}(1) \cdot p \lambda e_{\alpha} + p \lambda e_{-\alpha} - p \lambda e_{\alpha} + p^2X_{\mathcal{O}_F} \in L, \]
for any $\lambda \in \mathcal{O}_F$. Since $X_{\mathcal{O}_F}$ is spanned over $\mathcal{O}_F$ by $\{e_{\alpha}, h_{\alpha} \mid \alpha \in X \}$, it follows that $L$ contains $pX_{\mathcal{O}_F}/p^2X_{\mathcal{O}_F}$, as required. $\square$\\
\end{jpf}

\section{A normalisation lemma} \label{commutative algebra}

To finish the proof of \Cref{Main result canonical dimension one thm}, we need only establish \Cref{characteristic ideal fg element prop}. We will actually prove a stronger statement, which we may view as a relation of Noether's normalisation lemma.

\begin{jprop} \label{commutative algebra prop}
Let $K$ be any field, and $A$ be a commutative finitely-generated positively graded $K$-algebra of Krull dimension 1. Let $f \in A$ be a non-zero homogeneous element of positive degree and let $T = \{f^n \mid n \geq 0 \}$. Then $T^{-1}A$ is a finitely-generated module over the subalgebra $T^{-1}K[f]$.\\
\end{jprop}

When $G$ is a uniform pro-$p$ group with element $g \in G$, and $M$ is a $kG$-module of canonical dimension 1, we obtain \Cref{characteristic ideal fg element prop} by putting $A = \gr kG/J(M)$ and letting $f$ be the image of $\sigma(g-1)$ in $A$.\\

\begin{jpfof}[\Cref{commutative algebra prop}]
Consider the ideal of $f$-power torsion elements of $A$,
\[ I = \{ x \in A \mid f^nx =0 \text{ for some }n \geq 0\} \subseteq A, \]
then let $B=A/I$, and let $h=f+I \in B$. Localising the natural quotient map $A \rightarrow B$ induces an isomorphism of $T^{-1}A$-modules,
\[ T^{-1}A \cong T^{-1}B, \]
due to the definition of $I$. Hence, to show that $T^{-1}A$ is a finitely-generated $T^{-1}K[f]$-module, it is enough to show that $B$ is a finitely-generated $K[h]$-module.\\
By definition of $B$, the element $h \in B$ is not a zero divisor, so multiplication by $h^n$ induces an isomorphism of $B$-modules,
\[ B/Bh \cong Bh^n/Bh^{n+1}, \]
for all $n \geq 0$. Hence all of the subquotients in the descending sequence of $B$-modules
\[ B \geq Bh \geq Bh^2 \geq \dots \]
are isomorphic to $B/Bh$. It follows from the definition of Krull dimension that
\[\mathcal{K}(B/Bh) < \mathcal{K}(B) \leq \mathcal{K}(A)=1. \]
So $B/Bh$ is an Artinian $B$-module. Since $B$ is a finitely-generated commutative $K$-algebra, $B$ is Noetherian. Moreover, Hilbert's Nullstellensatz implies that all simple $B$-modules are finite-dimensional over $K$. It follows that $B/Bh$ is a finite-length $B$-module and finite-dimensional $K$-vector space.
Finally, it follows from the graded Nakayama lemma that $B$ is a finitely-generated $K[h]$-module, giving the result.\\
Namely, there is a finite-dimensional vector space $V \leq B$ such that 
\[ B = V + hB.\]
Let $B'=B/K[h]V$, a positively graded $K[h]$-module. Then $hB'=B'$, implying that for all $m \geq 0$, the graded pieces
\[ B'_m = (hB')_m = hB'_{m-\mathrm{deg}(f)},\]
so for all $n \geq 0$,
\[ B'_m = h^nB'_{m-\mathrm{deg}(f)n}. \]
Since $f$ is of positive degree and $B'$ is positively graded, it follows that $B'=0$. Thus $B = K[h]V$ is a finitely-generated $K[h]$-module, as required. $\square$\\
\end{jpfof}

I thank the reviewer for the elegant refinement of the statement and proof above.

\section{Krull dimension} \label{Krull dimension}

Following our result limiting the possible canonical dimensions of modules for $kG$, we can bound the Krull dimension of the Iwasawa algebra above. This is a straightforward deduction from a result for rings, \Cref{canonical dimension ring range bound cor}. We first state a fundamental characterisation of Krull dimension, see sections 6.1-6.2 of \jcite{mcconnell2001noncommutative}.

\begin{jlem} \label{Krull dimension fundamental lemma}
Let $R$ be a ring and $M$ be a left $R$-module. For any $n \geq 0$, the Krull dimension $\mathcal{K}(M) \leq n$ if and only if every strictly descending chain of submodules
\[ M=M_0 \geq M_1 \geq M_2 \geq \dots \]
satisfies $\mathcal{K}(M_j/M_{j+1}) \leq n-1$ for all but finitely many $j$.\\
\end{jlem}

If $M$ is Artinian and non-zero, then there are no strictly decreasing chains of submodules, so \Cref{Krull dimension fundamental lemma} implies $\mathcal{K}(M) = 0$ trivially. We can now give a simple result bounding Krull dimension by canonical dimension (recall the definition for general Auslander-Gorenstein rings from \jcite[Definition 4.5]{levasseur1992regulargraded}).

\begin{jlem} \label{Krull canonical dimension inequality lemma}
Let $R$ be an Auslander-Gorenstein ring and $M$ be a finitely-generated left $R$-module. Then $\mathcal{K}(M) \leq \Cdim(M)$.\\
\end{jlem}
\begin{jpf}
First note that Theorem 2.2 of \jcite{levasseur1992regulargraded} implies $\Cdim(M) = -\infty$ if and only if $M=0$. It follows from Proposition 4.5 of \jcite{levasseur1992regulargraded} that if $\Cdim(M)=0$, then $M$ is Artinian, meaning $\mathcal{K}(M)=0$. Thus we have proved the result in the case $\Cdim(M)=0$.\\
We now proceed by induction. Let $\Cdim(M)=n$ and let
\[ M = M_0 \geq M_1 \geq \dots \]
be a strictly descending chain of submodules. For all but finitely many $j$, $\Cdim(M_j/M_{j+1}) \leq n-1$, hence $\mathcal{K}(M_j/M_{j+1}) \leq n-1$, by Proposition 4.5 of \jcite{levasseur1992regulargraded} and induction. Then $\mathcal{K}(M) \leq n$ by \Cref{Krull dimension fundamental lemma}, as required. $\square$\\ 
\end{jpf}

We will not use the full strength of \Cref{Krull canonical dimension inequality lemma}, but include it since we believe it is of independent interest. In particular, \Cref{Krull canonical dimension inequality lemma} implies that the Krull dimension of $kG$ is at most $\dim G$ -- see also Theorem A of \jcite{ardakov2004krull}.\\
Although we could not find \Cref{Krull canonical dimension inequality lemma} in existing literature, we expect it is well known to experts. The next result improves the bound, at the cost of being more difficult to compute.

\begin{jprop} \label{canonical dimension module range bound prop}
Let $R$ be an Auslander-Gorenstein ring and $M$ be a Noetherian left $R$-module. Then 
\[ \mathcal{K}(M) \leq |\{\Cdim(N) \mid N\text{ is a non-zero subquotient of }M\}|-1. \]
\end{jprop}
\begin{jpf}
First let us define notation; for any Noetherian left $R$-module $X$, write
\[ S(X) = \{\Cdim(N) \mid N\text{ is a non-zero subquotient of }X\}, \]
and $s(X) = |S(X)|$. Now, if $\Cdim(M)=0$, we must have $S(M)=\{0\}$, and by \Cref{Krull canonical dimension inequality lemma}, 
\[ \mathcal{K}(M) \leq \Cdim(M) = 0 = s(M)-1, \]
as required.\\
We now proceed by induction on the canonical dimension. Let $\Cdim(M)=n$. By Proposition 4.5 of \jcite{levasseur1992regulargraded}, every strictly decreasing sequence of submodules
\[ M = M_0 \geq  M_1 \geq \dots \]
has $\Cdim(M_j/M_{j+1}) \leq n-1$ for all but finitely many $j$. Since all subquotients of $M$ are Noetherian, by induction
\[ \mathcal{K}\Big(\faktor{M_j}{M_{j+1}}\Big) \leq s\Big(\faktor{M_j}{M_{j+1}}\Big)-1, \]
for all such $j$. Moreover, $n \not \in S(M_j/M_{j+1})$, but
\[ S\Big(\faktor{M_j}{M_{j+1}}\Big) \cup \{n\} \subseteq S(M), \]
meaning that 
\[ \mathcal{K}\Big(\faktor{M_j}{M_{j+1}}\Big) \leq s\Big(\faktor{M_j}{M_{j+1}}\Big) -1 \leq s(M)-2\]
for all but finitely many $j$. Thus \Cref{Krull dimension fundamental lemma} implies
\[ \mathcal{K}(M) \leq s(M)-1, \]
giving the result. $\square$\\
\end{jpf}

For $R$ a ring, recall that $\mathrm{mod}(R)$ denotes the non-zero finitely-generated left $R$-modules.

\begin{jcor} \label{canonical dimension ring range bound cor}
Let $R$ be a Noetherian Auslander-Gorenstein ring. Then
\[ \mathcal{K}(R) \leq |\{\Cdim(M) \mid M \in \mathrm{mod}(R)\}|-1. \]
\end{jcor}
\begin{jpf}
Let $M$ be a finitely-generated (non-zero) $R$-module with generators $m_1, \dots, m_n$. Let $N = Rm_n$. Then by Proposition 4.5 of \jcite{levasseur1992regulargraded}, either $\Cdim(M) = \Cdim(N)$ or $\Cdim(M) = \Cdim(M/N)$. By induction on the number of generators, it follows that $M$ has the same canonical dimension as a cyclic $R$-module. Hence
\[ \{\Cdim(M) \mid M \in \mathrm{mod}(R)\} = \{\Cdim(M) \mid M\text{ is a non-zero subquotient of }R\}.\]
Then \Cref{canonical dimension module range bound prop} implies the result. $\square$\\
\end{jpf}

\Cref{Main result canonical dimension one thm} gives a restriction on the canonical dimension of modules for an Iwasawa algebra, hence \Cref{canonical dimension ring range bound cor} can be applied. We now deduce \Cref{Main result Krull dimension bound cor}.

\notinsubfile{\krulldimensionbound*}
\onlyinsubfile{
\begin{jcor} \label{Main result Krull dimension bound cor}
Let $\mathfrak{g}$ be a finite-dimensional $F$-Lie algebra which is semisimple and split over $F$. Let $G$ be a compact $p$-adic Lie group with associated Lie algebra $\mathcal{L}(G) \cong \mathfrak{g}$. If $[F:\mathbb{Q}_p] > 1$, the Krull dimension of $kG$ is at most $\dim G - 1$.\\
\end{jcor}
}
\begin{jpf}
If $M$ is a non-zero finitely-generated left $kG$-module, then $\Cdim(M)$ is an integer between 0 and $\dim G$. By \Cref{Main result canonical dimension one thm}, $\Cdim(M) \neq 1$, therefore by \Cref{canonical dimension ring range bound cor}, $\mathcal{K}(kG) \leq \dim G -1$. $\square$\\ 
\end{jpf}

We can combine this with Theorems A and B of \jcite{ardakov2004krull} to obtain the following.

\begin{jcor} \label{Krull dimension bound all fields cor}
Let $\mathfrak{g}$ be a finite-dimensional $F$-Lie algebra which is simple and split over $F$. Let $G$ be a compact $p$-adic Lie group with associated Lie algebra $\mathcal{L}(G) \cong \mathfrak{g}$. If $\mathfrak{g} \cong \mathfrak{sl}_2(\mathbb{Q}_p)$, then $\mathcal{K}(kG) = 3$. Otherwise,
\[ \lambda(\mathfrak{g}) \leq \mathcal{K}(kG) \leq \dim \mathfrak{g} -1, \]
where $\lambda(\mathfrak{g})$ is the maximum length of a chain of $\mathbb{Q}_p$-Lie subalgebras of $\mathfrak{g}$. \\
\end{jcor}

In fact the lower bound of \Cref{Krull dimension bound all fields cor} holds for any compact $p$-adic Lie group $G$, and Ardakov and Brown have conjectured that $\mathcal{K}(kG)=\lambda(\mathfrak{g})$. \Cref{Krull dimension bound all fields cor} means we can verify this in some cases. For example if $G=SL_3(\mathbb{Z}_p)$, then $\mathfrak{g} \cong \mathfrak{sl}_3(\mathbb{Q}_p)$, so $\dim_{\mathbb{Q}_p}(\mathfrak{g}) = 8$ and $\lambda(\mathfrak{g}) = 7$. It follows from \Cref{Krull dimension bound all fields cor} that $\mathcal{K}(kG) = 7$, see Corollary B of \jcite{ardakov2004krull}.\\

Our new result allows us to compute the Krull dimension of the Iwasawa algebra of a 6-dimensional group.

\begin{jcor} \label{quadratic extension Krull dimension cor}
Let $F$ be a quadratic extension of $\mathbb{Q}_p$ and $G=SL_2(\mathcal{O}_F)$. Then $\mathcal{K}(kG)=5$.\\
\end{jcor}
\begin{jpf} 
We have that $\mathfrak{g} = \mathcal{L}(G) \cong \mathfrak{sl}_2(F)$, and $\dim_{\mathbb{Q}_p} \mathfrak{g} = 6$, $\lambda(\mathfrak{g}) = 5$. By \Cref{Main result Krull dimension bound cor} and \cite[Theorem A]{ardakov2004krull}, it follows that $\mathcal{K}(kG)=5$. $\square$\\
\end{jpf}

Finally, we can also deduce a bound on the Krull dimension of the Iwasawa algebras of the compact general linear groups.

\begin{jprop} \label{Krull dimension GL_n equality prop}
Let $n \geq 2$ and $F$ be any finite field extension of $\mathbb{Q}_p$. The Krull dimension
\[ \mathcal{K}(kGL_n(\mathcal{O}_F)) = \mathcal{K}(kSL_n(\mathcal{O}_F)) + [F:\mathbb{Q}_p]. \]
\end{jprop}
\begin{jpf} 
Let $G = GL_n(\mathcal{O}_F)$. Let $I$ be the identity matrix, and 
\[ Z = \{ c I \mid c -1 \in p^{1+e}\mathcal{O}_F\} \leq G,\]
where $e=0$ if $p$ is odd, but $e=1$ if $p=2$. Then $Z$ is a closed central subgroup of $G$, and is uniform pro-$p$. We have
\[ \mathcal{L}(G) = \mathfrak{gl}_n(F), \quad \mathcal{L}(Z) = \{ c I \mid c \in F\}, \]
and there is a direct sum decomposition into Lie ideals,
\[ \mathfrak{gl}_n(F) = \mathfrak{sl}_n(F) \oplus \{ c I \mid c \in F\}. \]
Therefore
\[ \mathcal{L}(G/Z) = \mathcal{L}(G)/\mathcal{L}(Z) \cong \mathfrak{sl}_n(F). \]
Since $\mathcal{L}(SL_n(\mathcal{O}_F)) \cong \mathfrak{sl}_n(F)$, by \Cref{Lie algebra equivalence cor} it follows that 
\[ \mathcal{K}(k(G/Z)) = \mathcal{K}(kSL_n(\mathcal{O}_F)). \]
Recall the augmentation ideals 
\[ \epsilon_G(N) = kG\epsilon_N(N) = \epsilon_N(N)kG \]
for closed normal subgroups $N \leq G$, by Lemma 2.36 of \jcite{timmins2023coherence}. Now, $Z$ is a closed subgroup of the open uniform pro-$p$ normal subgroup
\[H = \{X \in G \mid X - I \in p^{1+e}\mathrm{Mat}_n(\mathcal{O}_F)\} \leq G.\]
In fact, the augmentation ideal $\epsilon_G(H)$ lies in the Jacobson radical of $kG$. To see this, if a maximal left ideal $M \subseteq kG$ does not contain $\epsilon_G(H)$, then
\[ M +  \epsilon_G(H) = M + \epsilon_H(H)kG = kG. \]
Since $kG$ is a finitely-generated $kH$-module and the Jacobson radical of $kH$ is $\epsilon_H(H)$, Nakayama's lemma implies $M=kG$, a contradiction. Thus $\epsilon_G(Z) \subseteq \epsilon_G(H)$ is contained in the Jacobson radical of $kG$. Now, $kZ$ is isomorphic to a commutative power series ring on $[F:\mathbb{Q}_p]$ elements, central in $kG$. Therefore the augmentation ideal $\epsilon_G(Z) = kG\epsilon_Z(Z)$ is generated by $[F:\mathbb{Q}_p]$ regular normal elements lying in the Jacobson radical of $kG$. By Corollary 1.9 of \jcite{walker1972localrings}, the Krull dimension of $kG$ is given by
\[ \mathcal{K}(kG) = \mathcal{K}(kG/\epsilon_G(Z)) + [F:\mathbb{Q}_p] = \mathcal{K}(k(G/Z)) + [F:\mathbb{Q}_p].\]
The result then follows since $\mathcal{K}(k(G/Z)) = \mathcal{K}(kSL_n(\mathcal{O}_F))$. $\square$\\
\end{jpf}

\begin{jthm} \label{Krull dimension bound general linear thm}
When $[F:\mathbb{Q}_p] > 1$ or $n >2$, the Krull dimension of $kGL_n(\mathcal{O}_F)$ satisfies
\[ \mathcal{K}(kGL_n(\mathcal{O}_F)) \leq n^2[F:\mathbb{Q}_p] -1. \]
Also, $\mathcal{K}(kGL_2(\mathbb{Z}_p))=4$, $\mathcal{K}(kGL_3(\mathbb{Z}_p))=8$, and $\mathcal{K}(kGL_2(\mathcal{O}_F))=7$ when $F$ is a quadratic extension of $\mathbb{Q}_p$.\\
\end{jthm}
\begin{jpf} 
By \Cref{Krull dimension bound all fields cor} and \Cref{Krull dimension GL_n equality prop}, if $n > 2$ or $F \neq \mathbb{Q}_p$,
\begin{align*} \mathcal{K}(kGL_n(\mathcal{O}_F)) &= \mathcal{K}(kSL_n(\mathcal{O}_F)) + [F:\mathbb{Q}_p]\\  &\leq \dim(\mathfrak{sl}_n(F)) -1 +  [F:\mathbb{Q}_p]. \end{align*}
But the $\mathbb{Q}_p$-vector space dimension $\dim(\mathfrak{sl}_n(F)) = (n^2-1)[F:\mathbb{Q}_p]$, so the claimed inequality follows.\\
The equalities in specific cases follow by combining \Cref{Krull dimension GL_n equality prop} with Corollary A of \jcite{ardakov2004krull}, Corollary B of \jcite{ardakov2004krull}, and \Cref{quadratic extension Krull dimension cor}, respectively. $\square$\\
\end{jpf}

\section{Holonomic representations} \label{Holonomic representations}

We now turn our attention to the mod $p$ representation theory of $p$-adic Lie groups. The theory of canonical dimension appears naturally when studying admissible representations. From now on, $G$ will denote an arbitrary $p$-adic Lie group unless otherwise stated. All our representations are over $k$, a field of characteristic $p$.\\

Recall that a smooth admissible representation $V$ of $G$ is one which satisfies
\[ V = \bigcup_{K \leq_o G } V^K ,\]
where $K$ ranges over all open subgroups of $G$, and every subspace of invariants
\[ V^K = \{ v \in V \mid x \cdot v = v\; \forall x \in K\}, \]
is finite-dimensional. Frequently we drop the adjective ``smooth''. For every admissible representation $V$, the $k$-linear dual
\[ V^\vee = \mathrm{Hom}_k(V,k) \]
is naturally a module for the \textit{augmented} Iwasawa algebra $kG$, which was studied in \jcite{kohlhaase2017smoothduality} and \jcite{timmins2023coherence}. Moreover, $V^\vee$ is a finitely-generated module for the Iwasawa algebra $kH$ of any compact open subgroup $H \leq G$. This property characterises the duals of admissible representations; by Corollary 1.8 and Proposition 1.9 of \jcite{kohlhaase2017smoothduality}, (see also \cite[2.2.12]{emerton2010ordinaryparts}), the dual map $V \rightarrow V^\vee$ induces an anti-equivalence of abelian categories
\[ \mathrm{Rep}_{G}^{\infty,a}(k) \longleftrightarrow \mathcal{M}_G, \]
where $\mathrm{Rep}_{G}^{\infty,a}(k)$ is the category of smooth admissible representations of $G$ over $k$, and $\mathcal{M}_G$ is as follows.

\begin{jdef} 
The category $\mathcal{M}_G$ is the full subcategory of $kG$-$\mathrm{Mod}$ whose objects are modules $M$ such that $M$ is finitely-generated as a $kH$-module for any, equivalently all, compact open subgroups $H \leq G$.\\
\end{jdef}

For any compact open subgroup $H \leq G$, any module $M$ in $\mathcal{M}_G$ is a finitely-generated $kH$-module, hence has a canonical dimension. The arguments of \Cref{change of group prop} and \Cref{Lie algebra equivalence cor} show that the canonical dimension is invariant with respect to the choice of $H$, thus we unambiguously denote it by $\Cdim(M)$.\\

We are now in a position to give the third characterisation of canonical dimension appearing in this article, which was originally proved in Proposition 2.18 of \jcite{emertonpaskunas2020density}.

\begin{jthm} \label{asymptotic growth thm}
Let $G$ be a $p$-adic Lie group and $V$ be a smooth admissible representation of $G$ over $k$. Let $K$ be an open uniform pro-$p$ subgroup of $G$.   There exists an integer $N \geq 0$ and real numbers $a,b > 0$ such that for all $n \geq N$,
\[ a(p^n)^d \leq \dim_k V^{K_n} \leq b (p^n)^d, \]
where $K_n = \{h^{p^n} \mid h \in K\}$ are the $p$-power subgroups of $K$, and the canonical dimension $\Cdim(V^\vee)=d$.\\
\end{jthm}

Thus, we see that the canonical dimension $\Cdim(V^\vee)$ associated to an admissible representation $V$ measures the growth rate of its invariant subspaces. In particular, $V$ is finite-dimensional if and only if $\Cdim(V^\vee)=0$.\\

When appropriate, we refer to $\Cdim(V^\vee)$ as the GK-dimension of $V$. We stress that this name is given only by analogy; \Cref{asymptotic growth thm} bears somewhat superficial similarity to the definitions of the usual Gelfand-Kirillov dimension. Iwasawa algebras are not Cohen-Macaulay with respect to the Gelfand-Kirillov dimension, and in all interesting cases the true Gelfand-Kirillov dimension of a $kH$-module will be infinite, as discussed in subsection 5.6 of \jcite{ardakov2006survey}. It is however the case that $\Cdim(V^\vee)$ is equal to the Gelfand-Kirillov dimension of the graded module $\gr V^\vee$, see Proposition 5.4 of \textit{loc.\! cit}.\\

\Cref{asymptotic growth thm} shows that the GK-dimension may be computed without reference to Iwasawa algebras, at least in principle. Of course, the results of this article would not have been obtained without this perspective, for example we have the following.

\begin{jthm} \label{non-compact canonical dimension one thm}
Let $\mathbb{G}$ be a split-semisimple affine group scheme defined over $F$, and $G=\mathbb{G}(F)$. Suppose $[F:\mathbb{Q}_p] > 1$. There are no admissible representations of $G$ over $k$ of GK-dimension 1.\\
\end{jthm}
\begin{jpf} 
A compact open subgroup of $G$ is $H = \mathbb{G}(\mathcal{O}_F)$. Then, $\mathcal{L}(H)$ is isomorphic to a finite-dimensional $F$-Lie algebra which is semisimple and split over $F$. By \Cref{Main result canonical dimension one thm}, $kH$ has no finitely-generated modules of canonical dimension 1. In particular, $\mathcal{M}_G$ has no objects of canonical dimension 1. Thus for every admissible representation $V$ of $G$, the GK-dimension $\Cdim(V^\vee) \neq 1$. $\square$\\
\end{jpf}

An immediate corollary to \Cref{non-compact canonical dimension one thm} is that any infinite-dimensional admissible representation of $G=\mathbb{G}(F)$ has GK-dimension at least two, when $F$ is a non-trivial extension of $\mathbb{Q}_p$. We expect that this lower bound will vary as $G$ varies, and in general will be higher than two.\\

We now turn our attention to representations that satisfy such a lower bound. Whenever a module for the Weyl algebra (or more generally, a $\mathcal{D}$-module) satisfies the lower bound of Bernstein's inequality, it is called holonomic, and holonomic modules are always of finite length - see \cite[Theorem 1.12, Definition 1.13, 1.16.1]{ehlers1987weyl}. We now prove the analogous result for representations of $p$-adic Lie groups.

\begin{jdef} 
Let $M$ be a module in $\mathcal{M}_G$. We say $M$ has minimal-positive dimension if $\Cdim(M) > 0$ and any other module $N$ in $\mathcal{M}_G$ has $\Cdim(N) \geq \Cdim(M)$ or $\Cdim(N)=0$.\\
\end{jdef}

\begin{jdef} 
Let $V$ be an admissible representation of $G$. $V$ is holonomic if and only if $V^\vee$ has minimal-positive dimension.\\
\end{jdef}

Notice that any admissible representation $V$ with GK-dimension $\Cdim(V^\vee)=1$ is automatically holonomic. We now prove the main result of this section. 

\notinsubfile{\holonomicmodulesfinitelength*}
\onlyinsubfile{
\begin{jthm} \label{Main result minimal-positive implies finite-length thm}
Suppose $G$ acts trivially on all of its finite-dimensional smooth representations over $k$. Every holonomic smooth admissible representation of $G$ is of finite length.\\
\end{jthm}
}
\begin{jpf}
Given a holonomic representation $V$, its dual $M=V^\vee$ has minimal-positive dimension. It is enough to show that $M$ has finite length, since the anti-equivalence of categories $\mathrm{Rep}_{G}^{\infty,a}(k) \leftrightarrow \mathcal{M}_G$ preserves simple objects.\\
Suppose $M$ is not of finite length. Now $M$ is Noetherian, because $M$ is a finitely-generated module over a Noetherian ring, hence $M$ must be non-Artinian. Let
\[ M=M_0 \geq M_1 \geq M_2 \geq \dots \]
be a strictly decreasing sequence of $kG$-submodules. By Proposition 4.5 of \jcite{levasseur1992regulargraded}, there exists $N \geq 0$ such that $\Cdim(M_n/M_{n+1}) < \Cdim(M)$ for all $n \geq N$. Since $M$ has minimal-positive dimension, $\Cdim(M_n/M_{n+1})=0$, thus $M_n/M_{n+1}$ is finite-dimensional for all $n \geq N$. Hence $M_N/M_n$ is finite-dimensional for all $n \geq N$.\\
Then $M_N=W^\vee$ for some admissible representation $W$ of $G$, with a strictly increasing sequence of finite-dimensional subrepresentations
\[ 0=W_N \leq W_{N+1} \leq W_{N+2} \leq \dots. \]
By assumption on $G$, $W_n \leq W^G$ for all $n \geq N$. But $W$ is admissible so $W^G$ itself is finite-dimensional, a contradiction. Thus $M$ is of finite length. $\square$\\
\end{jpf}

I thank the reviewer for the greatly shortened proof above. Some assumption on the group $G$ is necessary for \Cref{Main result minimal-positive implies finite-length thm} to hold, as illustrated by the example of the free module of rank one over the Iwasawa algebra $k\mathbb{Z}_p=k[[T]]$. In the next section we will show that Chevalley groups over $F$, for example $SL_n(F)$, satisfy the hypothesis of \Cref{Main result minimal-positive implies finite-length thm}, see \Cref{Chevalley fd module prop}.

\section{Finite length admissible representations} \label{Finite length GLn}

We can deduce interesting corollaries about certain admissible representations from \Cref{Main result minimal-positive implies finite-length thm}. First, we determine a large class of $p$-adic Lie groups whose finite-dimensional smooth representations are trivial -- we expect this is well-known to the experts.

\begin{jprop} \label{F fd module prop}
Let $F$ be a finite field extension of $\mathbb{Q}_p$ and $k$ be a field of characteristic $p$. All finite-dimensional smooth representations of the additive group $F$ are trivial.\\
\end{jprop}
\begin{jpf} 
Consider an $n$-dimensional smooth representation of $F$, $\phi: F \rightarrow GL_n(k)$, and let $a$ be a non-negative integer such that $n \leq p^a$. Since $\phi$ is a smooth representation, $\Ker \phi$ is an open subgroup of $F$.\\
Let $x \in F$. We will show that $\phi(x)$ equals the identity element $I$. Let $y = p^{-a}x \in F$. Since $\Ker \phi$ is an open subgroup, $p^{N}y \in \Ker \phi$ for some $N \geq 0$, and hence
\[ (\phi(y) -I)^{p^{N}} = \phi(y)^{p^{N}} -I = 0, \]
since $k$ is of characteristic $p$. Therefore $\phi(y)-I$ is a nilpotent matrix in $GL_n(k)$. By the Cayley-Hamilton theorem, it follows that $(\phi(y)-I)^n=0$, and hence as $n \leq p^a$,
\[ \phi(x) -I = (\phi(y)-I)^{p^a} = (\phi(y)-I)^n(\phi(y)-I)^{{p^a}-n} = 0, \]
so $\phi(x)=I$ as required. $\square$\\
\end{jpf}

I thank the reviewer for suggesting the improved proof of \Cref{F fd module prop} above.

\begin{jprop} \label{Chevalley fd module prop}
Let $X$ be a root system, $G=X(F)$ be the Chevalley group over $F$. All finite-dimensional smooth representations of $G$ are trivial. \\
\end{jprop}
\begin{jpf} 
By definition of the Chevalley group $X(F)$, $G$ is generated by subgroups
\[ \left(U_\alpha = \{x_\alpha(t) \mid t \in F\}\right)_{\alpha \in X},\]
and clearly $U_\alpha \cong F$ for all roots $\alpha \in X$. If $V$ is a smooth finite-dimensional representation of $G$ over $k$, then $V$ is a smooth finite-dimensional representation of each $U_{\alpha}$ under restriction. By \Cref{F fd module prop}, $V$ is a trivial representation for every $U_{\alpha}$. Since the $U_{\alpha}$ generate $G$, it follows that $V$ is a trivial representation of $G$. $\square$\\
\end{jpf}

The following application of \Cref{Main result minimal-positive implies finite-length thm} is immediate.

\begin{jcor} \label{X(F) finite length cor}
Let $X$ be a root system and $G=X(F)$ be the Chevalley group over $F$. Let $V$ be a smooth admissible representation of $G$ with $\Cdim(V^\vee)=1$. Then $V$ is of finite length.\\
\end{jcor}
\begin{jpf} 
Since $\Cdim(V^\vee)=1$, it is automatic that $V$ is holonomic. By \Cref{Chevalley fd module prop} and \Cref{Main result minimal-positive implies finite-length thm}, $V$ is of finite length. $\square$\\
\end{jpf}

If $X$ is indecomposable, notice that \Cref{Main result canonical dimension one thm} and \cite[Theorem 3.12, Proposition 3.15]{ardakov2004krull} tell us that no representations of GK-dimension one exist if $[F:\mathbb{Q}_p] > 1$ or $X \neq A_1$ -- that is, if $G \neq SL_2(\mathbb{Q}_p)$. See \Cref{Chevalley finite length thm} for a result on representations of GK-dimension two.\\ 

Frequently in the context of the (mod $p$) Langlands programme we are most interested in admissible representations of $GL_n(F)$. However, we do not possess any useful results limiting the canonical dimension of modules for $kGL_n(\mathcal{O}_F)$. In particular, every non-negative integer up to $[F:\mathbb{Q}_p]$ is realised as the canonical dimension of a $kGL_n(\mathcal{O}_F)$-module, by inflation along the homomorphism $kGL_n(\mathcal{O}_F) \rightarrow k\mathcal{O}_F^\times$ given by the determinant -- see \Cref{abelian example}. Moreover, $GL_n(F)$ has non-trivial finite-dimensional smooth representations, so does not satisfy the hypothesis of \Cref{Main result minimal-positive implies finite-length thm}.\\

Fortunately, many interesting representations of $GL_n(F)$ have a central character, and this allows us to apply Theorems \ref{Main result GL2(Qp) finite length cor} and \ref{Main result GLn finite length cor} by reducing to the case of special linear groups.

\begin{jdef} 
Let $V$ be a representation of a $p$-adic Lie group over the field $k$. $V$ has a central character if and only if every element of the centre $Z(G)$ acts on $V$ by a scalar in $k$.\\
\end{jdef}

\begin{jlem} \label{admissible restriction lemma}
Let $V$ be a smooth admissible representation of $GL_n(F)$ with a central character. Then $V$ is a smooth admissible representation of $SL_n(F)$ under restriction. Moreover $\Cdim_{GL_n(\mathcal{O}_F)}(V^\vee) = \Cdim_{SL_n(\mathcal{O}_F)}(V^\vee)$.\\
\end{jlem}
\begin{jpf} 
Let $Z$ be the centre of $GL_n(F)$. Consider the representation $V$ restricted to $Z$, giving a smooth representation $\chi: Z \rightarrow GL(V)$, which has open kernel because $Z$ acts by scalars. Let $K \leq Z$ be an open subgroup with $K \leq \Ker \chi$.\\
Let $H$ be an open subgroup of $SL_n(F)$. Since $K$ is open in $Z$, $HK$ is open in $SL_n(F)Z$, and hence in $GL_n(F)$. Because $V$ is admissible, $V^{HK}$ is finite-dimensional. But, $K$ acts trivially on $V$ by definition, so we have
\[ V^{HK} = (V^{K})^H = V^H. \]
Thus the fixed vectors of $V$ under any open subgroup of $SL_n(F)$ are finite-dimensional, so $V$ is an admissible representation of $SL_n(F)$.\\
Moreover, we can use the above equality of vector spaces to compare the canonical dimensions. Suppose that $H$ is uniform pro-$p$. Because $SL_n(F) \cap Z$ is finite, we may without loss of generality choose $K$ to be both uniform pro-$p$ and have trivial intersection with $H$. Then, $HK \leq GL_n(F)$ is an open uniform pro-$p$ subgroup, and recalling the notation of \Cref{asymptotic growth thm},
\[ (HK)_m = \{h^{p^m} \mid h \in HK \} = H_mK_m \]
for all $m \geq 0$. By the above, it follows that
\[ \dim_k V^{(HK)_m} = \dim_k V^{H_m} \]
for all $m \geq 0$. \Cref{asymptotic growth thm} shows that the canonical dimension of $V^\vee$ is determined by these integer-valued sequences of dimensions of invariant subspaces, hence the result follows. $\square$\\
\end{jpf}

\begin{jthm} \label{GL2(Qp) finite length cor}
Let $V$ be a smooth admissible representation of $GL_2(\mathbb{Q}_p)$, with central character, such that $\Cdim(V^\vee)=1$. Then $V$ is of finite length.\\
\end{jthm}
\begin{jpf} 
Since $V$ has a central character, \Cref{admissible restriction lemma} implies $V$ is a smooth admissible representation of $SL_2(\mathbb{Q}_p)$ under restriction, and 
\[ \Cdim_{SL_2(\mathbb{Q}_p)}(V^\vee)=\Cdim_{GL_2(\mathbb{Q}_p)}(V^\vee)=1. \]
By \Cref{X(F) finite length cor}, $V$ is a finite length representation of $SL_2(\mathbb{Q}_p)$, hence is of finite length as a representation of $GL_2(\mathbb{Q}_p)$. $\square$\\
\end{jpf}

The irreducible smooth admissible representations of $GL_2(\mathbb{Q}_p)$ were classified by Breuil in \jcite{breuil2003quelques}, and all have a central character by \jcite{berger2012centralcharacters}. The work of Morra in \jcite{morra2013invariant}, which builds on that of \Paskunas in \jcite{paskunas2010extensions}, shows that these representations correspond to modules of canonical dimension either one or zero. From this we can characterise the finite-length representations of $GL_2(\mathbb{Q}_p)$ with central character.

\notinsubfile{\GLQpfinitelength*}
\onlyinsubfile{
\begin{jcor} \label{classification finite length GL2(Qp) cor}
The finite length smooth admissible representations of $GL_2(\mathbb{Q}_p)$ with central character are precisely those with $\Cdim(V^\vee) \leq 1$.\\
\end{jcor}
}
\begin{jpf} 
By the results of Morra and Breuil, all irreducible admissible representations of $GL_2(\mathbb{Q}_p)$ have GK-dimension at most 1, since finite-dimensional representations have canonical dimension zero. Given any short exact sequence of admissible representations,
\[ 0 \rightarrow V_1 \rightarrow V_2 \rightarrow V_3 \rightarrow 0, \]
the equality $\Cdim(V_2^\vee) = \max(\Cdim(V_1^\vee), \Cdim(V_3^\vee))$ holds by Proposition 4.5 of \jcite{levasseur1992regulargraded}. Therefore all finite length admissible representations of $GL_2(\mathbb{Q}_p)$ have GK-dimension at most 1.\\
Conversely, if $V$ is a smooth admissible representation of $GL_2(\mathbb{Q}_p)$ of GK-dimension at most 1, then either $V$ is finite-dimensional, hence finite length, or $\Cdim(V^\vee)=1$, so $V$ is of finite length by \Cref{GL2(Qp) finite length cor}. $\square$\\
\end{jpf}

In \jcite{hu2024quaternion} many smooth admissible representations of (the units of) quaternion algebras over $\mathbb{Q}_p$ are exhibited, which have canonical dimension equal to one. As commented by the authors, these representations are not of finite length, which can occur since the quaternion algebras do not satisfy the condition of \Cref{Main result minimal-positive implies finite-length thm}.\\

We now use our new results on canonical dimension to find representations of finite length over groups of higher dimension.

\begin{jthm} \label{Chevalley finite length thm}
Let $X$ be an indecomposable root system, and $V$ be a smooth admissible representation of $X(F)$ where $[F:\mathbb{Q}_p] > 1$ or $X \neq A_1$. If $\Cdim(V^\vee)=2$, then $V$ is of finite length.\\
\end{jthm}
\begin{jpf} 
If $[F:\mathbb{Q}_p] > 1$, then by \Cref{non-compact canonical dimension one thm}, there are no admissible representations of $X(F)$ with GK-dimension $1$. Hence if $\Cdim(V^\vee)=2$, then $V$ is holonomic.\\
If $F=\mathbb{Q}_p$ and $X \neq A_1$, then by Theorem 3.12 and Proposition 3.15 of \jcite{ardakov2004krull}, combined with \Cref{Lie algebra equivalence cor}, the Iwasawa algebra $kX(\mathbb{Z}_p)$ has no modules of canonical dimension one. Since $X(\mathbb{Z}_p)$ is a compact open subgroup of $X(\mathbb{Q}_p)$, it follows that any admissible representation of $X(\mathbb{Q}_p)$ of GK-dimension 2 is holonomic.\\
In either case, \Cref{Chevalley fd module prop} and \Cref{Main result minimal-positive implies finite-length thm} then imply $V$ is of finite length. $\square$\\
\end{jpf}

In particular we may specialise \Cref{Chevalley finite length thm} to the case of $SL_n(F)$ -- we now apply this result.
\notinsubfile{\GLnfinitelength*}
\onlyinsubfile{
\begin{jcor} \label{GLn finite length cor}
Let $V$ be a smooth admissible representation of $GL_n(F)$ with central character. If $[F:\mathbb{Q}_p] > 1$ or $n >2$, and $\Cdim(V^\vee)=2$, then $V$ is of finite length.\\
\end{jcor}
}
\begin{jpf} 
By \Cref{admissible restriction lemma}, $V$ is an admissible representation of $SL_n(F)$ under restriction and $\Cdim_{SL_n(F)}(V^\vee) = \Cdim_{GL_n(F)}(V^\vee) = 2$. By \Cref{Chevalley finite length thm}, $V$ is a finite length representation of $SL_n(F)$, hence also of $GL_n(F)$. $\square$\\
\end{jpf}

Let us specialise to the case of the general linear group $GL_2(F)$. In \jcite{breuilXconjectures}, the following category of smooth admissible representations is considered.

\begin{jdef} 
Let $G=GL_2(F)$ where $F$ is a finite unramified extension of $\mathbb{Q}_p$. Consider the pro-$p$ Iwahori subgroup $I_1 \leq GL_2(F)$ and its centre $Z_1$. Let $J \trianglelefteq \gr k(I_1/Z_1)$ be the graded ideal given by equation (118) on page 114 of \jcite{breuilXconjectures}. We define the category $\mathcal{C}$ to be the full subcategory of $\mathrm{Rep}_{G}^{\infty,a}(k)$ whose objects $V$ have a central character, such that $\gr V^\vee$ is annihilated by some power of $J$.\\
\end{jdef}

Because $Z_1$ is pro-$p$, its only one-dimensional smooth representations are trivial, hence $V^\vee$ is a $k(I_1/Z_1)$-module, so the above is well-defined. Note that the category $\mathcal{C}$ considered above is not quite the same as the category with that name in \jcite{breuilXconjectures}. However, $\mathcal{C}$ is a full subcategory of that category, and the representations of most interest, namely those appearing in \cite[Theorem 1.1]{breuil2023gkdimension}, lie in $\mathcal{C}$.

\begin{jcor} \label{BHHMS rep cor}
Let $F$ be an unramified quadratic extension of $\mathbb{Q}_p$, and $G=GL_2(F)$. Let $V$ be a smooth admissible representation in the category $\mathcal{C}$. Then $V$ is of finite length.\\
\end{jcor}
\begin{jpf} 
The quotient ring $\gr k(I_1/Z_1)/J$ is a commutative algebra generated by $2[F:\mathbb{Q}_p]$ elements, with Krull dimension equal to $[F:\mathbb{Q}_p]$. Thus, when $H \leq I_1/Z_1$ is an open uniform pro-$p$ subgroup, $\gr kH/(\gr kH \cap J)$ is also a commutative algebra of Krull dimension $[F:\mathbb{Q}_p]$. Since $V$ is in $\mathcal{C}$, we have
\[ (\gr kH \cap J)^n \subseteq J^n \subseteq \mathrm{Ann}(\gr V^\vee) \]
for some $n \geq 1$. By \Cref{canonical dimension characteristic ideal thm},
\[\Cdim(V^\vee) \leq \mathcal{K}\Big(\faktor{\gr kH}{\gr kH \cap J}\Big)=[F:\mathbb{Q}_p]. \]
Since $[F:\mathbb{Q}_p]=2$, \Cref{non-compact canonical dimension one thm} and \Cref{admissible restriction lemma} imply that the GK-dimension $\Cdim(V^\vee)$ is either 0 or 2. By \Cref{Main result GLn finite length cor}, and since representations of GK-dimension 0 are finite-dimensional, it follows $V$ is of finite length. $\square$\\
\end{jpf}

It follows from \Cref{BHHMS rep cor} that the representations $\pi$ considered in \cite[Theorem 1.1]{breuil2023gkdimension} are of finite length when $F$ is a quadratic unramified extension of $\mathbb{Q}_p$.\\
This improves on some of the results of \jcite{breuilXconjectures}; when $\pi$ is as in subsection 3.3.2 of \textit{op.\! cit}, \Cref{BHHMS rep cor} implies that $\pi$ is of finite length, removing the condition $\pi^{K_1} \cong D_0(\bar{\rho})$ needed to prove this (in the case $\bar{\rho}$ reducible split) in subsection 3.3.5. That is, Theorem 1.3.8 of \textit{op.\! cit.\!} holds outside of the so-called ``minimal case'', if $F$ is a quadratic unramified extension of $\mathbb{Q}_p$.\\

Theorem 1.1 of \jcite{hu2022modpcohomology} states that certain representations $\pi(\bar{\rho})$ of $GL_2(F)$ (here $\bar{\rho}$ is non-semisimple) have GK-dimension $\Cdim(\pi(\bar{\rho})^\vee)=[F:\mathbb{Q}_p]$ when $F$ is an unramified extension of $\mathbb{Q}_p$. Again by \Cref{Main result GLn finite length cor}, if $F$ is a quadratic extension, then $\pi(\bar{\rho})$ is of finite length. Indeed, a stronger statement than this is shown in Theorem 1.7 of \textit{op.\! cit.\!} by different methods.\\

\bibliographystyle{cdim_jt.bst}
\bibliography{0_Whole_paper_Cdim.bib}

\end{document}